\documentclass[11pt]{article}
\usepackage{graphicx}
\usepackage{epstopdf}
\usepackage{amsmath, amssymb}
\usepackage{latexsym, color}
\usepackage{multirow}
\usepackage{mathtools}
\usepackage{xcolor}
\usepackage{hyperref}
\usepackage{graphicx}

\def\la{\big\langle}
\def\ra{\big\rangle}

\def\ds{\displaystyle}
\def\forall{\hbox{for all}~}
\def\L{{\bf L}}

\def\bfv{{\bf v}}

\def\bfn{{\bf n}}

\def\ve{\varepsilon}

\def\E{{\cal E}}
\def\A{{\cal A}}

\def\dint{\int\!\!\!\int}

\def\R{{\mathbb R}}

\def\F{{\cal F}}

\def\vs{\vskip 2em}

\def\v{\vskip 1em}
\def\O{{\cal O}}
\def\D{{\cal D}}

\def\C{{\cal C}}
\def\J{{\cal J}}
\def\M{{\cal M}}
\def\bega{\begin{array}}
\def\enda{\end{array}}
\def\begi{\begin{itemize}}
\def\endi{\end{itemize}}
\def\ov{\overline}
\def\curl{\hbox{curl }}

\def\Tilde{\widetilde}
\def\Hat{\widehat}

\def\avint{- \!\!\!\!\!\!\!\int}
\def\bpm{\begin{pmatrix}}
\def\epm{\end{pmatrix}}
\def\meas{\hbox{meas}}
\def\bel{\begin{equation}\label}
\def\eeq{\end{equation}}
\def\sqr#1#2{\vbox{\hrule height .#2pt
\hbox{\vrule width .#2pt height #1pt \kern #1pt
\vrule width .#2pt}\hrule height .#2pt }}
\def\square{\sqr74}
\def\endproof{\hphantom{MM}\hfill\llap{$\square$}\goodbreak}
\definecolor{cadmiumgreen}{rgb}{0.0, 0.42, 0.24}

\newtheorem{theorem}{Theorem}[section]

\newtheorem{proposition}{Proposition}[section]
\newtheorem{remark}{Remark}[section]
\newtheorem{definition}{Definition}[section]

\begin{document}

\title{\bf On the Optimal Control of Propagation Fronts}
\vs

\author{Alberto Bressan, Maria Teresa Chiri, and Najmeh Salehi\\
\,
\\
Department of Mathematics, Penn State University \\
University Park, Pa.~16802, USA.\\
\,
\\
e-mails: axb62@psu.edu, mxc6028@psu.edu, nfs5419@psu.edu.
}
\maketitle

\begin{abstract} We consider a controlled reaction-diffusion equation, motivated by 
a pest eradication problem.  Our goal is to derive a simpler model, describing 
the controlled evolution of a contaminated set. In this direction, the first part 
of the paper studies the optimal control of 1-dimensional traveling wave profiles.
Using Stokes' formula, explicit solutions  are obtained, which in some cases require
measure-valued optimal controls.
In the last section we introduce a family of optimization problems for a moving set.
We show how these can be derived from the original parabolic problems, by taking
a sharp interface limit.
\end{abstract}

\section{Introduction}
\label{sec:1}
\setcounter{equation}{0}

The control of parabolic equations is by now a classical subject
\cite{Coron, FM, LaTr, RBZ}.  More specifically, several studies have been devoted to 
the optimal harvesting
of spatially distributed populations \cite{CG, CGS, LM}.
Our present interest in the control of reaction-diffusion equations is
primarily motivated by models of pest eradication \cite{ACD, ACM,  CRo, SBGW}.
The controlled spreading of a population, in a simplest form, can be described 
by a semilinear parabolic equation
\bel{rd1}
u_t~=~f(u) + \Delta u  - g(u,\alpha).\eeq
Here  $u=u(t,x)$ denotes the 
population density at time $t$, at a location $x\in \R^2$.
The function $f$ describes the reproduction rate, while $\alpha=\alpha(t,x)$ is 
a distributed  control.  In a harvesting problem, the control function $\alpha$ 
accounts for the 
harvesting effort, while $g(u,\alpha)$ is the local amount of harvested biomass.
In the case of pest control, one may think of $\alpha(t,x)$ as the quantity of pesticides
sprayed at time $t$ at location $x$, while
$g(u,\alpha)$ describes the amount of population 
which is eliminated by this strategy.  We shall focus on the optimization problem
\begi
\item[{\bf (OP1)}] {\it Given an initial profile $u(0,x)=u_0(x)$ and a time interval $[0,T]$, determine a control $\alpha=\alpha(t,x)\geq 0$
so that, calling $u(t,x)$ the corresponding solution to (\ref{A6}), the total cost
\bel{rtcost} \J~\doteq~\int_0^T\phi\left(\int \alpha (t,x)\, dx\right)\, dt + \kappa_1\int_0^T \int u(t,x)\, dx\, dt + \kappa_2\,\int u(T,x)\, dx \eeq
is minimized.}
\endi
Here we think of $\int\alpha(t,x)dx$ as the {\it global control effort} at time $t$.

Several results are known on the existence of an optimal control, together with 
necessary conditions. However, one rarely finds explicit formulas, and 
optimal solutions can only be numerically computed. 
Aim the present paper is to derive a simplified model, for which optimal strategies
can be more easily found.  
By taking a sharp interface limit, our goal is to approximate the problem
{\bf (OP1)} with an optimal control problem for 
a moving set $\Omega(t)\subset \R^2$.  
Assuming that $f(1)=0$, $f'(1)<0$, so that 
$u=1$ is a stable equilibrium, we take
\bel{Omt}\Omega(t) ~=~\big\{ x\in \R^2\,;~u(t,x)\approx 1\big\}.\eeq
In connection with the cost functional (\ref{rtcost}),
in Section~\ref{sec:7} we will introduce 
a corresponding functional for the moving set $\Omega(t)$, and study its relation with 
{\bf (OP1)}. 

\begin{figure}[ht]
\centerline{\hbox{\includegraphics[width=13cm]{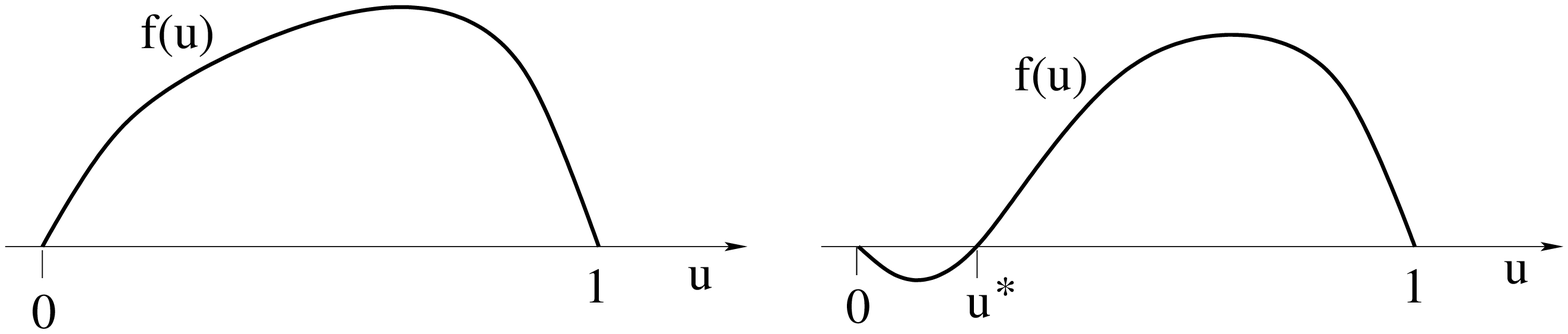}}}
\caption{\small Two possible shapes of the function $f$, satisfying the assumptions
{\bf (A1)} (monostable case),  and {\bf (A2)} (bistable case), respectively.}
\label{f:co1}
\end{figure}

Throughout the following, on the source function $f$ in (\ref{rd1}) we shall assume either one of the 
following  conditions (see Fig.~\ref{f:co1}):
\begi
\item[{\bf (A1)}] {\it $f\in \C^2$, and moreover
\bel{f1}
f(0)~=~f(1)~=~0,\qquad f''(u)\,<\,0\qquad\forall u\in [0,1].\eeq
}\endi
\begi
\item[{\bf (A2)}] {\it $f\in \C^2$, and moreover
\bel{f2}
f(0)~=~f(1)~=~0,\qquad f'(0)\,<\,0,\qquad f'(1)\,<\,0.\eeq
In addition, $f$ vanishes at only one intermediate point $u^*\in \,]0,1[\,$, where 
$f'(u^*)>0$.
}\endi
In addition, 
on the function $\phi$ we shall assume 
\begi
\item[{\bf (A3)}] {\it $\phi\in \C^2$, and moreover
\bel{pprop}
\phi(0)\,=\,0,\qquad \phi'(0)\geq 0,\qquad \phi''(s)>0\quad\forall s>0. \eeq
}\endi

Finally, we shall consider two simple choices of the function $g$ in (\ref{rd1}). Either
\bel{g1} g(u,\alpha) ~=~\alpha,\eeq
or else
\bel{g2} g(u,\alpha)~=~\alpha\, u.\eeq
In (\ref{g1}) the decrease of the pest population is proportional to the control effort.
On the other hand, (\ref{g2}) follows the more realistic harvesting  model
studied in \cite{BCS, CG, CGS}, where the local catch is proportional to the product
of the harvesting effort times the population density.

\begin{remark} {\rm As in \cite{BCS, CG, CGS}, the cost (\ref{rtcost}) has only linear growth 
w.r.t.~the local control effort $\alpha=\alpha(t,x)$.   Because of this, the optimal 
control may well be a measure,  not necessarily absolutely continuous w.r.t.~Lebesgue measure \cite{MR, R}.}
\end{remark}
In order to derive a cost functional for the motion of the set $\Omega(t)$ in (\ref{Omt}),
the key step lies in the analysis of traveling profiles for (\ref{rd1}).  Indeed,
the minimum cost associated to a traveling profile with speed $c$
will determine the local cost for 
moving the boundary $\partial\Omega(t)$ with speed $c$ in the normal direction.

The remainder of the paper is organized as follows.   Section~\ref{sec:2}
contains a brief review of the classical theory of traveling profiles for 1-dimensional reaction diffusion equations~\cite{F}.    In Section~\ref{sec:3} 
we consider traveling profiles having
a prescribed speed $c$  and requiring  a minimal control effort, i.e., minimizing the norm
$\|\alpha\|_{\L^1}$ of the control function in (\ref{rd1}).
We show that the above cost, associated with the traveling profile $u(t,x) = U(x-ct)$,
is computed by a line integral along the 
path $x\mapsto (U(x), U'(x))\subset\R^2$. Implementing 
a technique introduced in \cite{HH} (see also \cite{BoP}), one  can thus 
use Stokes' formula to estimate the difference
in cost between any two controlled  traveling profiles. In some cases, this allows us to 
explicitly determine the unique optimal profile.
In the remaining cases, in Section~\ref{sec:4} we prove the existence of a
 (possibly not unique)  optimal 
profile. Again, the optimal control $\alpha(\cdot)$ here can be a measure.

In Section~\ref{sec:5} we study how the minimum  cost $E(c)\doteq \min_\alpha \|\alpha\|_{\L^1}$ 
varies,
depending on the wave speed $c$.  As $c\to +\infty$, this cost always 
has linear growth. Indeed, an explicit formula  (\ref{EE})  for the
asymptotic behavior of the function $E$ can be given.  In the monostable case (\ref{f1}), 
we also show that this cost
is a convex function.  A partial extension of these results, to traveling profiles 
in a 2-dimensional space, is given in Section~\ref{sec:6}.

Section~\ref{sec:7} is the core of the paper.  Based on  the 
cost function $E(c)$  for optimal traveling profiles,
we introduce an optimization problem {\bf (OP2)} for
moving sets $t\mapsto\Omega(t)$.    
Our main result shows that this new cost (\ref{cost1}) can be attained as a limit
 of the costs corresponding to 
 a sequence of solutions of suitably rescaled parabolic problems.  

We remark that, in order to fully justify
 {\bf (OP2)} as a sharp interface limit of  {\bf (OP1)}, 
one should perform a detailed study of a corresponding  $\Gamma$-limit.
In the present paper,
the problem of characterizing the 
$\Gamma$-limit of the functionals in (\ref{Fen}) is left largely open.  
Under the assumptions {\bf (A1)}, two (small) steps in this direction are
worked out here. Proposition~\ref{p:53} proves the convexity of the
function $E(c)$. Moreover,  Proposition~\ref{p:81} shows that the 
optimal traveling profiles found in the 1-dimensional case are still optimal in two (or more) space dimensions.   Namely, more general traveling profiles of the form 
$u(x_1,x_2)~=~U(x_1-ct, x_2)$
do not achieve a lower cost, compared with profiles of the form
$u(x_1, x_2) = U(x_1-ct)$ which depend on the  single variable $x_1$.
For the definition and basic properties of $\Gamma$-limits we refer to \cite{Braides}.

Optimal control problems for moving sets, as in  {\bf (OP2)},
are studied in \cite{BCSa}, deriving necessary conditions for optimality and providing some
explicit formulas for the solutions.
Several other types of optimization problems for moving sets have been  considered in
\cite{Bblock, BMN, BZ, CPo, CLP}, motivated by different applications.

\section{Traveling wave solutions}
\label{sec:2}
\begin{figure}[ht]
\centerline{\hbox{\includegraphics[width=15cm]{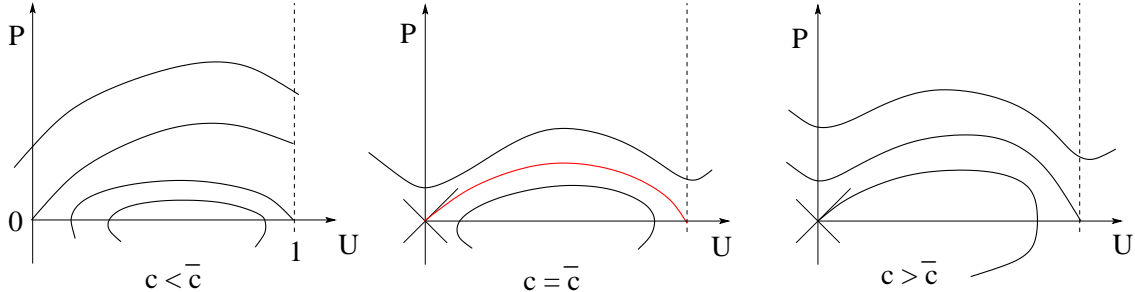}}}
\caption{\small A traveling profile for (\ref{A1}) corresponds to a heteroclinic orbit for the system
(\ref{T2}), connecting the points $(0,0)$ and $(1,0)$.  Under the assumptions {\bf (A2)}, 
such an orbit exists for one specific value $c=c^*$.}
\label{f:co5}
\end{figure}

\setcounter{equation}{0}
As a preliminary, we  review some basic facts on traveling waves for reaction-diffusion 
equations of the form
\bel{A1} u_t~=~f(u) + u_{xx}\,.\eeq
By definition, a traveling profile for (\ref{A1}) with speed $c$ is a solution of the form
\bel{TP} u(t,x)~=~U(x-ct).\eeq
This can be found by solving
\bel{TE} U''+cU' + f(U)~=~0.\eeq
Assuming that $f(0)=f(1)=0$,
we seek a solution $U:\R\mapsto [0,1]$ of (\ref{TE}) with asymptotic conditions
\bel{AC}
U(-\infty) ~=~0,\qquad U(+\infty) ~=~1.\eeq
Setting $P=U'$,  we thus need to find a heteroclinic orbit of the system
\bel{T2}
\left\{\bega{rl} U'&=~P,\\[1mm]
P'&=~-cP-f(U).\enda\right.\eeq
connecting the equilibrium points $(0,0)$ with $(0,1)$.
A phase plane analysis of the system (\ref{T2}) yields
\begin{theorem}\label{t:21} Consider the problem (\ref{TE})-(\ref{AC}).
\begi\item[(i)]   If $f$ satisfies {\bf (A1)}, then, for some number $c^* < 0$, 
there exists a traveling profile $U$ for every speed $c \leq c^*$.
\item[(ii)]  If $f$ satisfies {\bf (A2)}, then 
there exists a unique $c^*\in\R$ and a unique (up to a translation) traveling profile $U$
with speed $c=c^*$.
\endi
\end{theorem}
For a detailed proof, see Theorem~4.15 in \cite{F}.  In all cases, it can be shown that 
the traveling profile $U$ is monotone increasing.
A phase portrait of the system (\ref{T2}) in the bistable case  (\ref{f2})
is sketched in Fig.~\ref{f:co5}.

%
The Jacobian matrix at a point $(U,P)$ is
\bel{JM}J(U,P)~=~\begin{pmatrix}
0 && 1\\[1mm]
-f'(U) && -c\end{pmatrix}.\eeq
We observe that the assumption (\ref{f1}) implies that $(0,0)$ and $(1,0)$
are both saddle points in $(U,P)$ plane. The Jacobian matrix has real eigenvalues of opposite signs.   Indeed, solving
$$\lambda^2 + c\lambda + f'(U)~=~0$$
one obtains
\bel{la12}\lambda~=~{-c\pm\sqrt{c^2-4f'(U)}\over 2}\,.\eeq

We observe that, from (\ref{T2}), it follows
\bel{PE}{dP\over dU} ~=~-c - {f(U)\over P}\,,\qquad P(0)\,=\,P(1)\,=\,0.\eeq
Multiplying by $P$ and integrating over the interval $[0,1]$ one obtains
\bel{ws1}\int_0^1 P\, {dP\over dU}\, dU + \int_0^1 c P(U)\, dU~=~-\int_0^1 f(U)\, dU.\eeq
Therefore, the wave speed satisfies
\bel{ws2}c\, \int_0^1 P(U)\, dU~=~-\int_0^1 f(U)\, dU.\eeq
Since $U'=P>0$, this implies
\bel{ws3}\hbox{sign} ~c~=~-\hbox{sign} \int_0^1 f(U)\, dU.\eeq

\section{Optimal control of the wave speed}
\label{sec:3}
\setcounter{equation}{0}
In the setting of Theorem~\ref{t:21}, consider a speed $c>c^*$, so that the equation
(\ref{A1}) does not admit any traveling profile with speed $c$.
Given the function $g$ at (\ref{g1}) or (\ref{g2}), we then consider the controlled
equation
\bel{cwe}
u_t~=~f(u) + u_{xx} - g(u,\alpha).\eeq
Two questions now arise.
\begi
\item Does there exists a control $\alpha=\alpha(x-ct)\geq 0$ such that
(\ref{cwe})
admits a traveling wave solution with the prescribed speed $c$ ?
\item In the positive case, can this be achieved by an optimal control
$\alpha(\cdot)$, minimizing the
cost $\|\alpha\|_{\L^1}$ ?
\endi
Since the above cost has linear growth, to ensure the existence of an optimal  solution
we must reformulate the problem in a measure-valued setting.
Traveling profiles (\ref{TP}) thus correspond to solutions of 
\bel{tw3}U''+c U' + f(U)~=~\mu,\eeq
with asymptotic conditions
\bel{tw4} U(-\infty) \,=\,0,\quad U(+\infty) \,=\,1.\eeq

\begin{definition}
By $\M_c$ we denote the set of all 
positive, bounded Radon measures on the real line, such that
the problem (\ref{tw3})-(\ref{tw4}) has a monotonically increasing solution.
\end{definition}

More specifically, two optimization problems will be considered.
\begi
\item[{\bf (P1)}] {\it Assuming that  $f$ satisfies either {\bf (A1)} or {\bf (A2)}, 
given  a speed $c> c^*$,
find a measure $\mu\in \M_c$ which minimizes }
\bel{J0} J_0(\mu)\,\doteq\,\mu(\R).\eeq

\item[{\bf (P2)}] {\it In the bi-stable case where $f$ satisfies {\bf (A2)}, 
given a speed $c> c^*$,
find a measure $\mu\in \M_c$ which minimizes }
\bel{J1} J_1(\mu)\,\doteq\, \int_{\R} {1\over u}\, d\mu.\eeq
\endi

\subsection{The optimal solution for problem (P1).}
Setting $P=U'$,  a solution to  (\ref{tw3})-(\ref{tw4}) corresponds
to a solution of  \bel{T7}
\left\{\bega{rl} U'&=~P,\\[1mm]
P'&=~-cP-f(U) +\mu,\enda\right.\eeq
starting at  $(0,0)$ and reaching $(1,0)$.
Since $f$  is bounded and the measure $\mu$ has finite total mass, 
any such solution will have bounded total variation.

We observe that, at a point $\bar x$ where $\mu$ concentrates a positive mass, 
the derivative $P$ has an upward jump:
$$U'(\bar x+) - U'(\bar x-)~=~\mu(\{\bar x\})~>~0.$$
Following \cite{BR,R}, to the graph
$$\bigl\{(U,P)=(U(x), U'(x))\,;~~x\in\R\big\}$$
we add a (finite or countable) set of
 vertical segments,
at places where $P=U'$ has an upward jump.    By a suitable parameterization, 
this yields a Lipschitz curve 
\bel{gams} s~\mapsto ~\gamma(s)~=~\bigl(U(s), P(s)\bigr),
\qquad\quad s\in [0, \bar s],\eeq containing the graph of 
 the solution of (\ref{T7}).

  The cost $J_0$ in (\ref{J0}) can now be expressed as 
\bel{J01}
J_0(\gamma)~=~\int_0^{\bar s} \left[\Big( {f(U)\over P } + c\Big) U'(s) +P'(s)  \right]ds
~=~\int_\gamma\left[ \Big( {f(U)\over P } + c\Big)dU+dP \right].\eeq
This is to be minimized over a family $\A_c$ of admissible curves, defined as follows.

\begin{definition}
Given a wave speed $c$, we call $\A_c$ the set of all 1-Lipschitz curves
of the form $s\mapsto \gamma(s) = (U(s), P(s)\in \R^2$ such that, for some interval $[0, \bar s]$, one has
\bel{ad1} \gamma(0)=(0,0),\qquad \gamma(\bar s) = (1,0),
\qquad
P(s)\geq 0\qquad \forall  s\in [0,\bar s],\eeq
\bel{ad2}  \bigl|\gamma(s_1)-\gamma(s_2)\bigr|~\leq~ |s_1-s_2| \qquad\forall  s_1,s_2\in [0,\bar s],
\eeq
\bel{ad3} U'(s)\,\geq\,0,
\qquad P'(s)~\geq~\bigl( - f(U(s)) - c\bigr)U'(s)\qquad\hbox{for a.e.~}s\in [0, \bar s].\eeq
\end{definition}

\begin{figure}[ht]
\centerline{\hbox{\includegraphics[width=8cm]{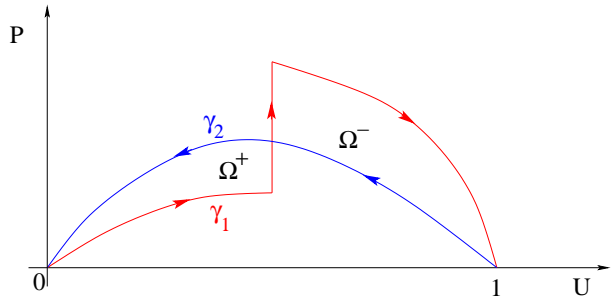}}}
\caption{\small  Estimating the difference in cost (\ref{Stokes}) between 
two paths $\gamma_1$ and $\gamma_2$. }
\label{f:sc16}
\end{figure}

Following an idea introduced in \cite{HH}, we use Stokes' theorem to compute the difference in cost between any two paths $\gamma_1,\gamma_2\in \A_c$.  Defining the vector field
$$\bfv~=~\left( {f(U)\over P}+c\,,~1\right),$$
by (\ref{J01}) we obtain
\bel{Stokes} J_0(\gamma_1) - J_0(\gamma_2)~=~
\left(\int_{\gamma_1}  - \int_{\gamma_2} \right)\bfv~=~\left( \dint_{\Omega^+} -
\dint_{\Omega^-} \right)\omega.\eeq
Here 
\bel{curl} \omega~= ~\curl \bfv~=~{f(U)\over P^2}\,,\eeq
while $\Omega= \Omega^+ \cup \Omega^-$ is the region enclosed between the two curves.
As shown in Fig.~\ref{f:sc16}, we call $\Omega^+$ the portion of this region whose boundary is traversed counterclockwise,
and $\Omega^-$ the portion whose boundary is traversed clockwise, when traveling first along $\gamma_1$, then along $\gamma_2$.

The formula (\ref{Stokes}) allows to immediately determine the optimal traveling 
wave profile, for the cost functional $J_0$.  

\begin{figure}[ht]
\centerline{\hbox{\includegraphics[width=17cm]{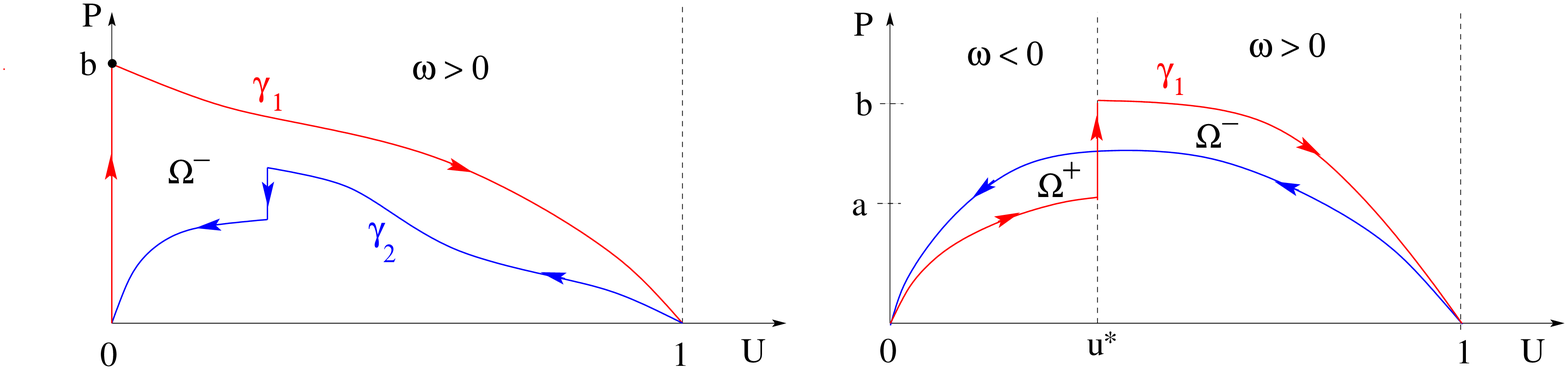}}}
\caption{\small  Left: if $f$ is always positive, then the curve $\gamma_1$
yields a lower cost, compared with any other admissible curve $\gamma_2$.    
Right: if $f$ is negative for $u<u^*$ and positive for $u>u^*$, then 
the curve $\gamma_1$ is the optimal one. }
\label{f:co16}
\end{figure}

\begi
\item[(i)] In the monostable case, where $f$ satisfies (\ref{f1}), we have
$\omega\geq 0$ throughout the domain.
As shown in Fig.~\ref{f:co16}, left, 
consider the solution of (\ref{T2}) through the saddle point $(1,0)$.
Since we are assuming $c>c^*$, this solution will cross the $P$-axis at some point 
$b>0$.   We then take $\gamma_1\in \A_c$ to be the path consisting of 
a vertical segment from $(0,0)$ to $(0,b)$, together with the trajectory of (\ref{T2})
from $(0,b)$ to $(1,0)$.   We claim that $\gamma_1$ is optimal.

Indeed, let $\gamma_2\in \A_c$ be any other admissible path.
Our definition of $\gamma_1$ implies that, by moving first along $\gamma_1$ then 
along $\gamma_2$, the boundary of the region enclosed by the two curves
is traversed clockwise.   Hence (\ref{Stokes}) implies
\bel{Stok2} J_0(\gamma_1) - J_0(\gamma_2)~=~-
\dint_{\Omega^-} \omega~=~-\dint_{\Omega^-} {f(U)\over P^2} \, dPdU~\leq~0.\eeq
Hence $\gamma_1$ is optimal.

\item[(ii)] In the bistable case, where $f$ satisfies (\ref{f2}), the function $\omega$ is
negative for $u<u^*$ and positive for $u>u^*$.
As shown in Fig.~\ref{f:co16}, right, 
let $(u^*, a)$ be the point reached by the trajectory of (\ref{T2}) 
through $(0,0)$, when it crosses the vertical line $\{U=u^*\}$.    Similarly, let $(u^*, b)$ be the point reached by the trajectory of (\ref{T2}) 
through $(1,0)$, when it crosses the vertical line $\{U=u^*\}$.  
Since we are assuming $c>c^*$, it follows that $a<b$.

Define $\gamma_1\in \A_c$ to be the path obtained by concatenating these two trajectories,
together with a vertical segment joining $(u^*,a)$ with $(u^*,b)$.
 We claim that $\gamma_1$ is optimal.

Indeed, let $\gamma_2$ be any other admissible path, connecting $(0,0)$ with $(1,0)$.
Our definition of $\gamma_1$ implies that, by moving first along $\gamma_1$ then 
along $\gamma_2$, the boundary of the region enclosed by the two curves
is traversed counterclockwise for $u<u^*$ and clockwise for $u>u^*$.   Hence (\ref{Stokes}) implies
\bel{Stok3} J_0(\gamma_1) - J_0(\gamma_2)~=~\left( \dint_{\Omega^+}-\dint_{\Omega^-} \right)
 \omega~=~\left( \dint_{\Omega^+}-\dint_{\Omega^-} \right){f(U)\over P^2} \, 
dPdU~\leq~0,\eeq
because the ratio $f(U)/P$ is negative on $\Omega^+$ and positive on $\Omega^-$.
Hence $\gamma_1$ is optimal.
\endi

We summarize the above analysis, stating the results in the original coordinates $u=u(t,x)$.
Consider the problem of minimizing the total mass $\mu(\R)$ among all positive measures
for which the equations (\ref{tw3})-(\ref{tw4}) have a solution.

\begin{theorem}\label{t:31}
For every $c>c^*$, the  problem {\bf (P1)} has a unique solution (up to translations).  
\begi
\item[(i)] In the monostable case, where $f$ satisfies (\ref{f1}), the optimal 
traveling profile can be uniquely determined by the equations
\bel{Uopt1} \left\{
\bega{rll}
U(x)&=~0\qquad &\hbox{if}~~x\leq 0,\\[2mm]
U'' + cU' + f(U)&= ~0\qquad &\hbox{if}~~x\leq 0,\\[2mm]
\ds   \lim_{x\to +\infty} U(x)&=~1.\enda
\right.\eeq

\item[(ii)]  In the bistable case, where $f$ satisfies (\ref{f2}), the optimal 
traveling profile can be uniquely determined by the equations
\bel{Uopt2} \left\{
\bega{rl} U'' + cU' + f(U)&= ~0\qquad\hbox{separately for $x<0$ and for $x>0$},\\[2mm]
U(0)&=~u^*, \\[2mm]
 \lim_{x\to -\infty} U(x)&=~0,\qquad   \lim_{x\to +\infty} U(x)~=~1\,. \enda
\right.\eeq
\endi
In both cases one has $f(u(0))=0$, and the optimal measure 
$\mu$ is a point mass located at the origin.  The minimum cost is
\bel{minco}
C_{min}~=~\mu\bigl(\{0\}\bigr)~=~ U'(0+)- U'(0-).\eeq
\end{theorem}

\begin{figure}[ht]
\centerline{\hbox{\includegraphics[width=15cm]{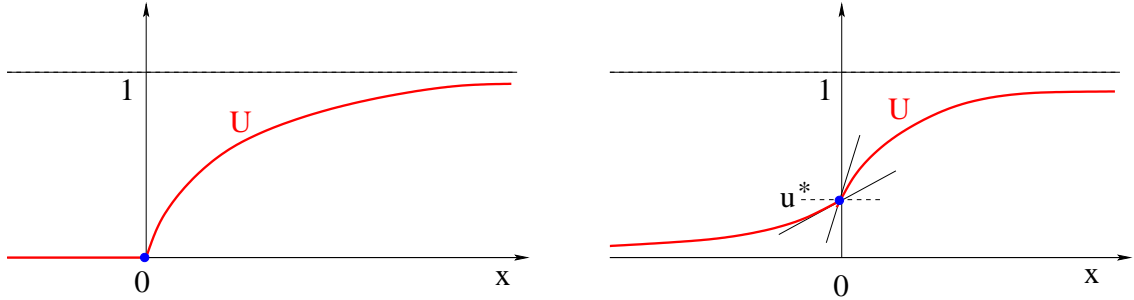}}}
\caption{\small The optimal traveling profiles described in Theorem~\ref{t:31}.
Left: the profile (\ref{Uopt1}). Right: the profile (\ref{Uopt2}).
}
\label{f:sc17}
\end{figure}

\subsection{The optimal solution for problem (P2).}
Next, consider the bistable case, but with cost functional (\ref{J1}).
Integrating along paths in the $U$-$P$ plane, instead of (\ref{J01}) we now find
\bel{J11}
J_1(\gamma)~=~\int_0^{\bar s} {1\over U} \left[\Big( {f(U)\over P } + c\Big) U'(s) +P'(s)  \right]ds
~=~\int_\gamma\left[ \Big( {f(U)\over U  P } + {c\over U} \Big)dU+{1\over U} dP \right].\eeq
Again, this is to be minimized among all admissible curves $\gamma\in \A_c$.
 Defining the vector field
$$\bfv~=~\left( {f(U)\over U P}+c\,,~{1\over U}\right),$$
and recalling (\ref{J11}) we now obtain
\bel{Stok4} J_1(\gamma_1) - J_1(\gamma_2)~=~
\left(\int_{\gamma_1}  - \int_{\gamma_2} \right)\bfv~=~\left( \dint_{\Omega^+} -
\dint_{\Omega^-} \right)\omega,\eeq
where now
\bel{curl2} \omega~= ~\curl \bfv~=~{f(U)\over U P^2}- {1\over U^2}  \,.\eeq
The region where $\omega>0$ is found to be
\bel{posr2}{\cal D}^+~\doteq~
\Big\{ (U,P)\,;~\omega(U,P)>0\Big\}~=~\Big\{ (U,P)\,;~0<P< P^*(U) \Big\}.\eeq
where 
\bel{P*} P^*(U)~\doteq~\sqrt{U\, f(U)}.\eeq

\begin{figure}[ht]
\centerline{\hbox{\includegraphics[width=10cm]{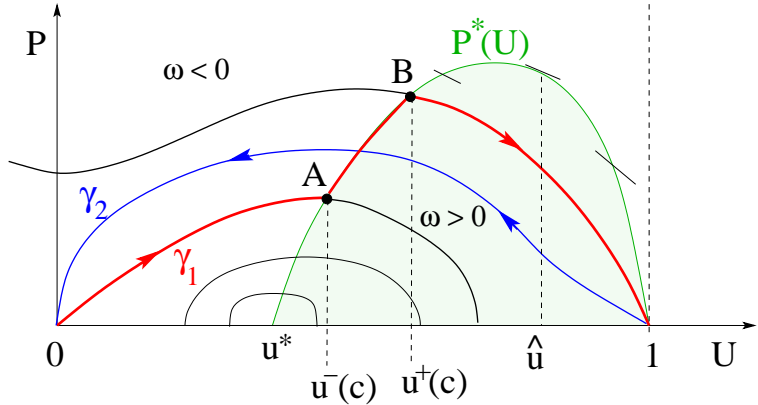}}}
\caption{\small  By Stokes' theorem, the path $\gamma_1$ going from $(0,0)$ to
$A$, then from $A$ to $B$, then from $B$ to $(1,0)$ has a lower cost than 
any other admissible path $\gamma_2\in \A_c$.}
\label{f:co17}
\end{figure}

Consider the situation shown in Fig.~\ref{f:co17}.    
Let $\gamma_1$ be the path obtained by concatenating:
\begi
\item The  trajectory of (\ref{T2}) exiting from $(0,0)$, until it reaches a point $A$
on the curve where $P=P^*(U)$.
\item  The  trajectory of (\ref{T2}) starting from  $(1,0)$, and moving backwards
until it reaches a point $B$
on the curve where $P=P^*(U)$.
\item The arc of the curve where $P=P^*(U)$, between $A$ and $B$.
\endi

Assume that the abode two trajectories of (\ref{T2}), passing through the points 
$A$ and $B$ respectively, do not have further
intersections with the  curve $P=P^*(U)$, for $u^*<U<1$.
Then $\gamma_1$ is optimal.


%

%

We give here a sufficient condition that guarantees that every trajectory
of (\ref{T2}) can cross the curve $P=P^*(U)$ only twice, thus ruling out
the configuration in Fig.~\ref{f:co19}.

Along this curve we have
\bel{dp*}{d\over dU} {P^*(U)}~=~{f(U) + U f'(U)\over 2\sqrt{ U f(U)}}\,.\eeq
Writing the equation (\ref{T7}) in the form
\bel{T77}
P'~=~-c - {f(U)\over P} + z^*(U),\eeq
a direct computation shows that, in the region where $P=P^*$ we must have
\bel{zopt} z^*(U)~=~{3 f(U) + U f'(U)\over 2\sqrt{U f(U)}} + c\,.\eeq
This leads us to consider the function
\bel{defg}2g(u)~\doteq~\bigl[3f(u)+ u f'(u)\bigr]\cdot \bigl[uf(u)\bigr]^{-1/2},\eeq
and seek a condition that will ensure that this function is monotonically decreasing.
A straightforward differentiation yields
$$\bega{l} 2 g'(u)~=~\bigl[3f'(u)+  f'(u)+ u f''(u)\bigr]\cdot \bigl[uf(u)\bigr]^{-1/2}\\[3mm]
\qquad \qquad\qquad \ds-{1\over 2} \bigl(3f(u)+ u f'(u)\bigr)\bigl( f(u) + u f'(u)\bigr)\cdot \bigl[uf(u)\bigr]^{-3/2}
\\[3mm]
\ds~ =~\bigg[ \Big(4u f(u)f'(u)+ u^2 f(u) f''(u)\Big) - {1\over 2} \Big( 3 f^2(u) + 4 u f(u) f'(u) + u^2 (f'(u))^2\Big)
\bigg] \cdot \bigl[uf(u)\bigr]^{-3/2}\\[4mm]
\ds~=~\bigg[ 2u f(u)f'(u)+ u^2 f(u) f''(u)-  {3\over 2} f^2(u) -{1\over 2} u^2 (f'(u))^2
\bigg] \cdot \bigl[uf(u)\bigr]^{-3/2}.
\enda
$$
Hence the inequality we need is
\bel{fneq}
-3 f^2(u) - u^2 (f'(u))^2+ 4u f(u)f'(u)+ 2u^2 f(u) f''(u)~\leq~0.\eeq

By the inequality $2\sqrt{3}uf(u)f'(u)\leq 3f^2(u)+u^2(f'(u))^2$, it follows that (\ref{fneq}) holds if, in addition to {\bf (A2)}, the function $f$ satisfies:
\begi
\item[{\bf (A4)} ] For all $u\in [u^*, 1]$, one has $ (4-2\sqrt 3) f'(u) + 2uf''(u)\leq 0$.
 \endi

\begin{theorem}\label{t:32} 
	Let $f$ be a function satisfying the assumptions {\bf (A2)} and {\bf (A4)} 
	Then the inequality (\ref{fneq}) holds for every $u\in [0,1]$.
	
Moreover, for every wave speed $c>c^*$, the optimization problem {\bf (P2)} admits a unique solution.
	The optimal measure
	is absolutely continuous w.r.t.~Lebesgue measure. There exists
	two points $u^*<u^-(c)<u^+(c) <1$ such that   the optimal solution in 
	(\ref{T77}) has the form
	\bel{occ}z^*(u)~=~\left\{\bega{cl} \ds{3 f(U) + U f'(U)\over 2\sqrt{U f(U)}} + c
	\qquad&\hbox{if} \quad u\in [u^-(c),u^+(c)],\\[3mm]
	0\qquad&\hbox{otherwise.}\enda\right.\eeq
\end{theorem}

{\bf Proof.}
As shown in Fig.~\ref{f:co17}, in the region
$$J^+~\doteq~\bigl\{ u\in [u^*,1]\,;~z^*(u)>0\bigr\}$$
where the control is strictly positive, the graph of the function $P^*$ intersects
transversally 
all trajectories of the system (\ref{T2}).   Therefore, if we can prove that the set $J^+$
is an interval, say $J^+ = \,]u^*, \Hat u[\,$ as shown in Fig.~\ref{f:co17},  we are done.
We start observing that, under the assumptions {\bf (A2)} on $f$, the function $z^*$ defined 
at (\ref{zopt}) satisfies
$$z^\ast (u^\ast)>0 \quad \hbox{ and }\quad  z^\ast (1)<0, $$
hence by continuity of $z^\ast$ there exists at least one point $\bar u\in [u^\ast, 1]$ such that $z^\ast(\bar u)=0$. We claim  that for any $c>c^*$ this point is unique. Indeed, consider the function 
\bel{ef} h(u)~\doteq ~- {3\,f(u)+u\,f'(u)\over 2\, \sqrt{u\,f(u)}}\,,\qquad  \quad u\in [u^\ast, 1].\eeq 
Since
\bel{lg}
\lim_{u\rightarrow u^*+}h(u)~=\,-\infty, \qquad \qquad \lim_{u\rightarrow 1-}h(u)~=\,+\infty \, ,
\eeq 
our claim will be proved by showing that 
$h$ is strictly increasing. Indeed, this is true because $h=-g$, with $g$ defined in (\ref{defg}).




\begin{remark} {\rm There is a large set of functions satisfying {\bf (A2)} and {\bf (A4)}.
For example,  the cubic $f(u)=-x(x-1)(x-2/3)$  satisfies {\bf (A4)} with a strict inequality. Therefore the same is true for any small perturbation  in $\C^2_0\bigl([0,1]\bigr)$.
}
\end{remark}

\v
\section{Existence of optimal traveling profiles}
\label{sec:4}
\setcounter{equation}{0}

\begin{figure}[ht]
\centerline{\hbox{\includegraphics[width=12cm]{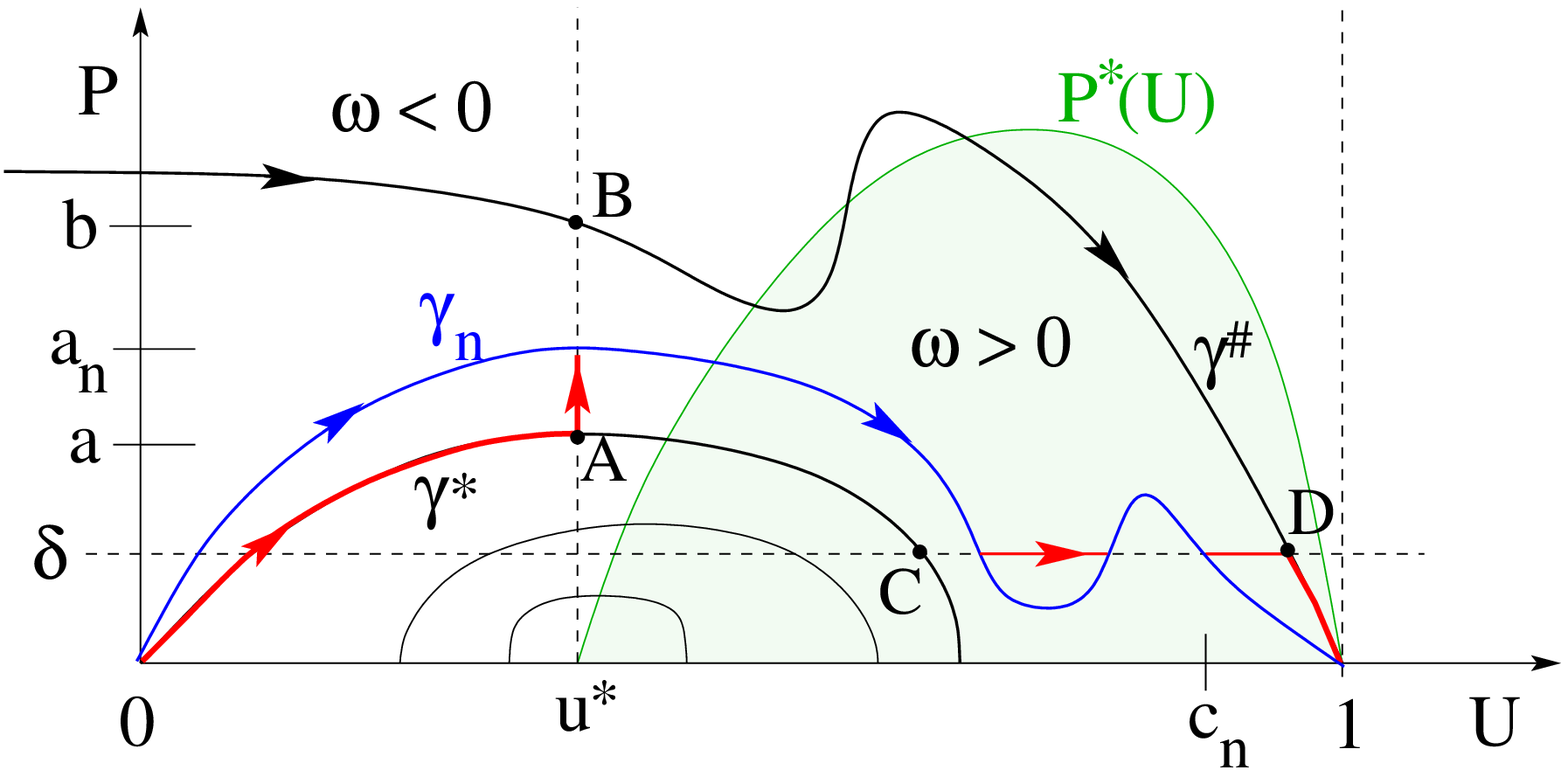}}}
\caption{\small The construction in step 2 of the proof of Theorem~\ref{t:32}.
An admissible path $\gamma_n$ is replaced by a path having smaller cost, and having
the additional properties (i)--(iii).
Notice that this is a case where the trajectory $\gamma^\sharp$ of (\ref{T2}) through the point
$(1,0)$ has multiple intersections with the curve $P=P^*(U)$.  When this happens,
Theorem~\ref{t:32} cannot be applied.}
\label{f:co19}
\end{figure}

For more general source functions $f$, satisfying {\bf (A2)} but not {\bf (A4)}
we still prove existence of an optimal measure $\mu$ yielding the traveling profile.
However, in the situation shown in Fig.~\ref{f:co19}
the structure of this measure can be more complicated than in
the case covered by  Theorem~\ref{t:31}.

\begin{theorem}\label{t:41}
Let $f$ satisfy the assumptions {\bf (A2)}, and let $c^*$ be the speed of a traveling wave for 
(\ref{A1}). Then, for every $c\geq c^*$ the minimization problem {\bf (P2)} has a measure valued solution.
\end{theorem}

{\bf Proof.} {\bf 1.}  
As shown in Fig.~\ref{f:co19}, let $\gamma^*$ be the trajectory of (\ref{T2})
originating from $(0,0)$, and let $\gamma^\sharp$ the trajectory of (\ref{T2})
reaching $(1,0)$.    Notice that every admissible path $\gamma\in \A_c$
is contained in the region bounded by $\gamma^*$, $\gamma^\sharp$, and the
$U$-axis.

Under the assumption $c\geq c^*$, a path with finite cost does exist.
Indeed, let $A=(u^*,a)$ and $B=(u^*,b)$ be the points where
the trajectories $\gamma^*, \gamma^\sharp$ cross  the vertical line $\{U=u^*\}$,
respectively.    Then the path $\gamma$ obtained by concatenating
\begi\item the portion of $\gamma^*$ from $(0,0)$ to $A$,
\item the vertical segment from $A$ to $B$, and
\item the portion of $\gamma^\sharp$ from $B$ to $(1,0)$
\endi
is an admissible path with cost
 $J_1(\gamma) = {b-a\over u^*}<+\infty$.   Notice that this is the path that minimizes the 
cost functional $J_0$, but of course it may not be optimal for $J_1$.
\v
{\bf 2.}
We can now consider a minimizing sequence of
paths $\gamma_n\in \A_c$, say $$s~\mapsto~ \gamma_n(s)~=~ \bigl( U_n(s), P_n(s)\bigr),
\qquad 
s\in [0, s_n],\quad  n\geq 1,$$
such that 
$$\lim_{n\to\infty} ~J_1(\gamma_n)~=~\inf_{\gamma\in\A_c} J_1(\gamma).$$

Adapting the arguments used in the previous section, based on Stokes' theorem, 
we now replace each path $\gamma_n$ with a modified path $\Tilde \gamma_n$
having some additional properties. 
As shown in Fig.~\ref{f:co19}, 
let $(u^*,a_n)$  be the first point where
the path $\gamma_n$ intersects the vertical line $U=u^*$.    We can then replace the 
portion of $\gamma_n$ between $(0,0)$ and $(u^*, a)$ with the portion 
of $\gamma^*$ from $(0,0)$ to the point $A=(u^*,a)$, 
together with a vertical segment joining $(u^*,a)$ with $(u^*, a_n)$.

Next, consider the portion of $\gamma_n$ in a neighborhood of the terminal point $(1,0)$.
Since the measure $\mu$ is positive,  that this portion must lie below the trajectory $\gamma^\sharp$ of (\ref{T2}) 
through $(1,0)$.
Moreover, in a neighborhood of $(1,0)$ the path $\gamma^\sharp$ lies below the curve $P=P^*(U)$. Indeed, in view of (\ref{dp*}),
$$\lim_{U\to 1-} {d\over dU} P^*(U)~=~\lim_{U\to 1-} \left( {\sqrt {f(U)}\over 2U} + {\sqrt U \, f'(U)\over 2 f(U)}\right)~=~-\infty,
$$
while, along $\gamma^\sharp$, by (\ref{JM})-(\ref{la12}) we have
$$\lim_{U\to 1-} {dP\over dU}~=~{-c-\sqrt{ c^2 - 4 f'(1)}\over 2}\,.$$

We now choose $\delta>0$ small enough so that, calling
$C,D$ the points where the horizontal line $\{P=\delta\}$ intersects
the trajectories $\gamma^*$ and $\gamma^\sharp$ respectively,
one has
$$-f(U)   - cP ~\leq~0$$
along the horizontal segment with endpoints $C,D$.    In other words, 
all trajectories of (\ref{T2}) cross this segment downward.   

Call $(c_n, \delta)$ the last point where the path $\gamma_n$ crosses the horizontal line
$\{P=\delta\}$. 
We then replace the last portion of $\gamma_n$ with a horizontal segment 
joining $(c_n,\delta)$ with  $D$, together with the 
portion of trajectory $\gamma^\sharp$ joining $D$ with $(1,0)$.
Furthermore, we replace any additional portions of the path $\gamma_n$
lying below the line $\{P=\delta\}$ with horizontal segments.

After these modifications, we obtain a new path $\Tilde\gamma_n$.
Since the function $\omega$ at (\ref{curl2}) is negative on the strip where $U<u^*$,
by (\ref{Stok3}) we have
$$J_1(\Tilde \gamma_n)\,\leq\, J_1(\gamma_n).$$
In view of the above construction we can now assume that
every path $\gamma_n$ in our minimizing sequence has the following properties:
\begi\item[(i)] The initial portion of $\gamma_n$  coincides with the path $\gamma^*$, from
$(0,0)$ to the point $A=(u^*,a)$.
\item[(ii)] The final portion of $\gamma_n$ coincides with the path $\gamma^\sharp$, from
the point $D$ to $(1,0)$. 
\item[(iii)] The intermediate portion of $\gamma_n$, 
between $A$ and $D$, remains inside the 
domain where $U\in [u^*,1]$ and $P\geq\delta$.
\endi
\v
{\bf 3.} 
  By parameterizing each path $\gamma_n$ by arc-length, we can assume that all 
maps $\gamma_n$ are 1-Lipschitz and that
the intervals $[0, s_n]$ are uniformly bounded.  
By possibly taking a subsequence, and using Ascoli's theorem, we achieve
the convergence
$$s_n\to \bar s, \qquad \gamma_n(s)\to \gamma(s)\qquad \forall s\in [0, \bar s[\,.$$
Moreover,  for any fixed $\ve>0$ the convergence is uniform on the subinterval where
$s\in [0, \bar s-\ve]$.
\v
{\bf 4.} We claim that the limit path is admissible, namely $\gamma\in\A_c$.
Indeed, the identities (\ref{ad1}) are clear.
Moreover, the limit of 1-Lipschitz curves is still 1-Lipschitz, hence (\ref{ad2})
holds as well.   Finally, we observe that the  differential constraint (\ref{ad3}) can be
formulated in terms of the  differential inclusion
\bel{DIG}\gamma' (s)~\in~F(\gamma(s)),\eeq
where
$$F(U,P)~=~\Big\{ (\dot u,\dot p)\in \R^2\,;~\dot u^2+\dot p^2\leq 1\,,\quad \dot u\geq 0,~~
\dot p ~\geq~
(- f(U) - c) \dot u\Big\}.$$
Since the multifunction $F$ is continuous, with compact, convex values, 
the set of solutions to the differential inclusion (\ref{DIG}) is closed under uniform convergence~\cite{AC}.  This 
shows that $\gamma\in \A_c$.
\v
{\bf 4.} It now remains to show that the limit path $\gamma$ is optimal.
Namely
\bel{J33}
J_1(\gamma)~=~\int_\gamma\left[ \Big( {f(U)\over U  P } + {c\over U} \Big)dU+{1\over U} dP \right]~=~\lim_{n\to\infty} J_1(\gamma_n).\eeq
Using (\ref{Stok3})-(\ref{curl2}) we obtain
\bel{omin} \lim_{n\to\infty} \Big( J_1(\gamma)- J_1(\gamma_n)\Big)~=~\lim_{n\to\infty}
\left( \dint_{\Omega^+_n} -
\dint_{\Omega^-_n} \right)\left({f(U)\over U P^2}- {1\over U^2}\right)~=~0.\eeq
Indeed, by construction, for every $n\geq 1$ the region $\Omega_n=\Omega_n^+ \cup
\Omega_n^-$ enclosed between the two curves is contained within
the region where $U\geq u^*$ and $P\geq \delta$.   On this region, the 
integrand in (\ref{omin})  is continuous and uniformly bounded.   Since the area of $\Omega_n$
shrinks to zero, we conclude that the above limit vanishes, proving the optimality of    
$\gamma$.\endproof

The above theorem provides the existence of an optimal 
profile, but it does not guarantee its uniqueness (up to translation).  From step {\bf 2}
of the proof, we can obtain some information about the optimal measure $\mu$.
Namely $\mu$ is supported on a region where $U\in [u^*, 1-\ve]$, for some $\ve>0$.
In particular, the optimal profile coincides with a solution of (\ref{TE})
for $U\in \,]0, u^*]$ and for $U\in [1-\ve, 1[\,$.

\section{The minimum cost, depending on the wave speed}
\label{sec:5}
\setcounter{equation}{0}
In setting considered in Theorem~\ref{t:32}, the optimal control has the form
(\ref{occ}).   Calling $[u^-(c), u^+(c)]$ the interval where
the control is nonzero (see Fig.~\ref{f:co17}), the minimum cost is thus
\bel{minc}\bega{rl}
E(c)&\ds=~\int_{u^-(c)}^{u^+(c)} {z^*(u)\over u }\, du
~=~\int_{u^-(c)}^{u^+(c)} {3 f(u) + u f'(u) + 2c \sqrt{u f(u)}\over 2u\sqrt{u f(u)}}\, du
\\[4mm]&\ds=~{3\over 2} \int_{u^-(c)}^{u^+(c)} {\sqrt{f(u)}\over u\sqrt{u}}\, du +{1\over 2} 
\int_{u^-(c)}^{u^+(c)}  { f'(u)\over \sqrt{uf(u)}}\, \, du + c \ln\left({u^+(c)\over u^-(c)}\right).
\enda\eeq
The first two terms on the right hand side of (\ref{minc}) are uniformly bounded.   The last term is the only one that grows without bound, as $c\to +\infty$.  The next proposition yields
more precise information on the asymptotic behavior of $E(c)$.
\v
\begin{proposition}\label{p:51} Let $f$ be a function satisfying the assumptions {\bf (A2)} and {\bf (A4)}.
As in Theorem~\ref{t:32},  let $z^*$  in (\ref{occ}) be the optimal control for the problem 
 {\bf (P2)}.
Then, as $c\to +\infty$, one has $u^-(c)\to u^*$, $u^+(c)\to 1$.
Moreover, the function
$E(c)$ in (\ref{minc}) has the asymptotic behavior
\bel{EE}
E(c)~=~ c\, | \ln u^* |+\int_{u^*}^1 \left({3\sqrt{f(u)}\over2 u\sqrt{u}} + { f'(u)\over 2 \sqrt{uf(u)}}\right)  du  + e(c), \eeq
where the additional term has size $e(c)=\O(1) \cdot {1\over c}$.
\end{proposition}

{\bf Proof.} All of the above conclusions will be proved by showing that 
there exists a constant $\alpha>0$ such that
\bel{upmc}u^-(c)-u^*~\leq~{\alpha\over c^2}\,,\qquad\qquad 1-u^-(c)~\leq~{\alpha\over c^2}\,.\eeq
\v
{\bf 1.}
Consider the equation 
\bel{PE1}{dP\over dU} ~=~-c - {f(U)\over P},\,\eeq
associated with the system (\ref{T2}).
One can observe that for $c>0$, 
\bel{SD} {dP\over dU}~\leq~ -c \qquad \hbox{ on } \text{  }[u^\ast, 1] .\eeq
Therefore $P$ attains its maximum at some point in the interval  $[0, u^\ast]$, where the right hand side of (\ref{PE1}) vanishes.   
For $u \in [0, u^*]$ one has $f_{\text{min}} \le f(u) \le 0$. Therefore
\bel{Umax}0~\leq~ P(U)~\leq~ -{f_{\text{min}}\over c}\qquad\qquad\forall U\in [0,1]. \eeq
In turn, (\ref{Umax}) implies
\bel{u-}-{f_{min}\over c}~\geq~ \sqrt{{u^-(c)f( u^-(c))}}.
\eeq
As $c\to +\infty$, both sides of (\ref{u-}) approach zero, hence $f(u^-(c))\to 0$ and 
$u^-(c)\to u^*$.   Moreover, performing a Taylor approximation at $u=u^*$, for a suitable constant $\alpha>0$ we find
\bel{u*-}\sqrt{ u^-(c)-u^*} ~\leq~{\sqrt{\alpha}\over c}\,.\eeq
proving the first inequality in (\ref{upmc}).

\v
{\bf 2.} 
To achieve an estimate on $u^+(c)$  we observe that, for every $c>0$, the function
$$Z(U)~=~c (1-U) \qquad\qquad U\in [u^*, 1]\,,$$
is a subsolution of
$${d\over dU} P(U)~=~ - c - {f(U)\over P}\,,\qquad\qquad P(1)=0.$$
Indeed, $$Z'~=~-c~\geq~- c - {f(U)\over Z} \qquad\qquad U\in [u^*, 1]\,.$$

In turn, this implies
$$P^*(u^+(c))~=~\sqrt{ u^+(c)) f(u^+(c))}~\geq~Z(u^+(c))~=~ c\bigl(1-u^+(c)\bigr).$$
Once again, for a suitable constant $\alpha>0$ this implies
$$\sqrt {\alpha (1-u^+(c))}~\geq~c\bigl(1-u^+(c)\bigr),$$
proving the second inequality in (\ref{upmc}).
\v
{\bf 3.} The asymptotic expansion (\ref{EE}) is now a consequence of 
(\ref{minc}), together with the inequalities in (\ref{upmc}).
\endproof
\v

\begin{figure}[ht]
\centerline{\hbox{\includegraphics[width=10cm]{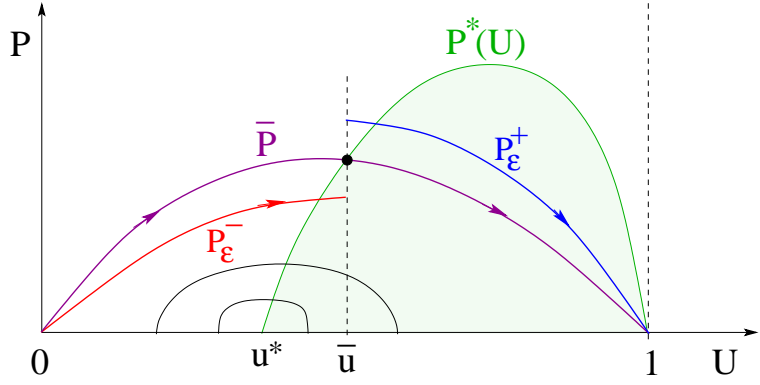}}}
\caption{\small When $c=c^*$ there exists a 
heteroclinic orbit $P=\ov P(U)$ of (\ref{T2}) through $(0,0)$ and $(1,0)$.   When $c=c^*+\ve$ with $\ve>0$,
one obtains an unstable manifold $P=P^-_\ve(U)$ through $(0,0)$, and  a 
stable manifold $P=P^+_\ve(U)$ through $(1,0)$.}
\label{f:co28}
\end{figure}

Next, we analyze the behavior of $E(c)$ as $c\downarrow c^*$.
Going back to the system (\ref{T2}), we observe that the unstable manifold
through $(0,0)$ and the stable manifold through $(1,0)$ depend continuously on the parameter $c$.   This impliies
\bel{lupm}
\lim_{c\to c^*} u^-(c)~=~
\lim_{c\to c^*} u^+(c)~=~\bar u.\eeq
As shown in Fig.~\ref{f:co28}, here $\bar u$ is the point where
 the heteroclicic orbit $P=\ov P(U)$
connecting  $(0,0)$ with $(1,0)$  intersects the graph of the function 
$P=P^*(U)\doteq\sqrt{Uf(U)}$,
in the case $c=c^*$.   

Setting $c=c^*+\ve$, we  now denote by $P=P^-_\ve(U)$ and $P=P^+_\ve(U)$ the 
corresponding unstable and stable manifolds  through $(0,0)$ and   $(1,0)$, respectively
(see Fig.~\ref{f:co28}).  By definition, the functions $P^-_\ve, P^+_\ve$ thus provide the solutions to 
\bel{dpu}{dP\over dU} ~=~-{ f(U)\over P}- (c^*+\ve),\eeq
respectively with boundary conditions
$$P_\ve^-(0)~=~0,\qquad\qquad P_\ve^+(1)~=~0.$$
Differentiating (\ref{dpu}) w.r.t.~the parameter $\ve$, we obtain
the asymptotic expansions
\bel{aex}
P^-_\ve(U)~=~\ov P(U) +\ve Y^-(U) + o(\ve),\qquad\qquad 
P^+_\ve(U)~=~\ov P(U) +\ve Y^+(U) + o(\ve),\eeq
where $o(\ve)$ denotes a higher order term, as $\ve\to 0$.
Here the functions $Y^-,Y^+$ are determined  by solving the linearized equations
\bel{Ylin} {dY\over dU}~=~{f(U)\over [\ov P(U)]^2} \, Y- 1,\eeq
with boundary conditions
$$\lim_{U\to 0+} Y^-(U)~=~0,\qquad\qquad \lim_{U\to 1-} Y^+(U)~=~0,$$
respectively.  In view of the formula (\ref{minc}), we now obtain

\begin{proposition}\label{p:52} In the same setting as Proposition~\ref{p:51}, 
we have the asymptotic expansion
\bel{Ec0}E(c^*+\ve)~=~\ve\cdot {Y^+(\bar u) - Y^-(\bar u)\over \bar u} + o(\ve).\eeq
where $o(e)$ denotes a higher order infinitesimal as $\ve\downarrow 0$.
\end{proposition}
{\bf Proof.}  
For notational convenience, set
$$v^-~\doteq~{d\over d\ve} u^-(c^*+\ve)\bigg|_{\ve=0}\,,\qquad\qquad v^+~\doteq~{d\over d\ve} u^+(c^*+\ve)\bigg|_{\ve=0}\,.$$
Differentiating w.r.t.~$\ve$ the identities
$$P_\ve^-\bigl(u^-(c^*+\ve)\bigr)~=~P^*\bigl(u^-(c^*+\ve)\bigr),\qquad\qquad
P_\ve^+\bigl(u^-(c^*+\ve)\bigr)~=~P^*\bigl(u^+(c^*+\ve)\bigr),$$
and using (\ref{aex}), we obtain
\bel{Ypm}Y^\pm(\bar u) + \left( -c^* - {f(\bar u)\over\bar u}\right) v^\pm~=~(P^*)'(\bar u)\cdot v^\pm.\eeq
It is now convenient to write the minimum cost (\ref{minc}) in the form
\bel{mc2}
E(c)~=~\int_{u^-(c)}^{u^+(c)} {z^*(u)\over u}\, du~=~
\int_{u^-(c)}^{u^+(c)} {1\over u}\left( c+ {f(u)\over P^*(u)} + (P^*)'(u)\right)\, du\,.
\eeq
Differentiating (\ref{mc2}) 
w.r.t.~$c$, when  $c=c^*$ and $u^-(c^*)= u^+(c^*)=\bar u$ we obtain
$$E'(c^*)~ =~ {1\over \bar u}\left( c^*+ {f(\bar u)\over P^*(\bar u)} + (P^*)'(\bar u)\right)\cdot (v^+-v^-)~=~{ Y^+(\bar u) -
Y^-(\bar u)\over \bar u}\,,
$$
where the second identity follows from (\ref{Ypm}).  This yields (\ref{Ec0}).
\endproof
\v
The last result in this section is concerned with minimum cost  $J_0$ in (\ref{J0}), 
as a  function of the wave speed $c$,
but now in the mono-stable case (\ref{f1}).
In this case, we can prove that the function $E(c)$ is convex.

\begin{proposition}\label{p:53} Consider the minimization problem {\bf (P1)}, assuming that
 $f$ satisfies {\bf (A1)}.   Then the minimum cost $E(c)$ is an increasing, convex function
 of the speed $c\in [c^*, +\infty[\,$. 
\end{proposition}

{\bf Proof.} As shown by Theorem~\ref{t:31},
in this case the optimal control consists of a point mass at the origin.
The minimum cost is thus simply
$P(0)$, where $P=P(U)$ is the solution to
\bel{PU3} {dP\over dU}~=~-c P - {f(U)\over P},\qquad\qquad P(1)=0.\eeq
We shall write $P=P(U,c)$ to emphasize the dependence of the solution on the additional
parameter $c$.   Our main concern is the convexity of the map
$c\mapsto P(0,c)$.
To understand this issue, we set~
$w = {\partial P\over\partial c}$.  
This function satisfies the linear ODE
\bel{wode}
{dw\over dU}~=~-1 + {f(U)\over P^2(U)}\, w,\qquad\qquad w(1)=0.\eeq
If now $c_1<c_2$, then $P(U, c_1)<P(U, c_2)$. Therefore
$$-1 + {f(U)\over P^2(U, c_1)} ~>~-1 + {f(U)\over P^2(U, c_2)}\,.$$
In view of (\ref{wode}), this yields
\bel{cvex}
w(U, c_1)~<~w(U, c_2)\qquad \qquad \forall U\in [0, 1[\,.\eeq
showing that the map $c\mapsto P(U,c)$ is convex, for every $U\in [0,1]$.
\endproof

%

\section{Traveling profiles in two space dimensions}
\label{sec:6}
\setcounter{equation}{0}

In Theorem~\ref{t:31} we proved that, for any speed $c>c^*$, the optimization
problem {\bf (P1)} admits a unique optimal solution. The optimal measure is a point mass
located at a point where $f(u)=0$.  

Aim of this section is to prove a similar result for traveling waves in two space
dimensions.   
We thus consider the corresponding parabolic equation on the 2-dimensional strip
$\bigl\{ (x_1,x_2)\in \R\times \,]0,1[\,\bigr\}$, namely
\bel{A9} u_t~=~f(u) + \Delta u- z(x),\eeq
with Neumann boundary conditions:
\bel{NBC}
u_{x_2}(x_1,0)~=~u_{x_2}(x_1,1)~=~0\qquad\forall  x_1\in \R.\eeq
Given a speed $c>c^*$, we consider a traveling wave profile  $u=u(x_1, x_2)$
which satisfies 
\bel{2dtv}
f(u) + c u_{x_1} + \Delta u - z~=~0,\eeq
together with (\ref{NBC}) and with limits 
\bel{l12}  \lim_{x_1\to -\infty} u(x_1,x_2)~=~0,\qquad \lim_{x_1\to +\infty} u(x_1,x_2)~=~1.
\eeq
Among all such profiles, 
obtained by different choices of the function $z=z(x_1, x_2)\geq 0$,
we seek  to
minimize  the total effort
\bel{m2}\|z\|_{\L^1} ~\doteq~\int_{\R\times [0,1]} \bigl| z(x)\bigr|\, dx.  \eeq
We claim  that, even by choosing control functions $z$ which depend 
on both variables $x_1, x_2$, one cannot achieve a smaller cost compared with
the 1-dimensional case, where $z$ is a function of the variable $x_1$ alone.

\begin{proposition}\label{p:81} Let $f$ satisfy the assumptions
{\bf (A1)}. Given $c>c^*$, 
let $u=u(x_1, x_2)$ be a solution to  (\ref{NBC})--(\ref{l12}), for some
 nonnegative smooth function $z=z(x_1,x_2)$.
Calling $C_{min}$ the minimum cost  for the 1-dimensional problem at (\ref{minco}),
one has
\bel{lbcost} \|z\|_{\L^1}~\geq~C_{min}\,.\eeq
\end{proposition}
\v
{\bf Proof.}
{\bf 1.} By (\ref{2dtv}) and the boundary conditions (\ref{NBC}), (\ref{l12}),
one has
\bel{mc}\|z\|_{\L^1}~=~\int_{\R\times [0,1]}  f(u)\, dx + c.\eeq
Introducing  the level sets
$$\Sigma(s)~\doteq~\Big\{ x\in \R\times [0,1]\,;~u(x)=s\Big\},$$
the integral on the right hand side of (\ref{mc}) can be written as 
\bel{eff2}\int_0^1 \left( \int_{\Sigma(s)} {1\over |\nabla u(x)|} \, d\ell
\right) f(s)\,ds.\eeq
Here $d\ell$ denotes the arc-length along the level curve $\Sigma(s)$.
%
\v
{\bf 2.} Assuming that $f$ satisfies (\ref{f1}),
let
$U:\R_+\mapsto [0,1]$ be the 
optimal traveling profile constructed at (\ref{Uopt1}).
For every $s\in [0,1]$, define
\bel{pP1}\psi(s)~\doteq~\int_{\{u(x)>s\}}f(u(x))\, dx,\qquad\qquad \Psi(s)~\doteq~\int_{\{U(x)>s\}}f(U(x))\, dx\,.
\eeq
We claim that, for every $s\in [0,1]$,
\bel{comp1}
\psi(s)~\geq~\Psi(s).\eeq
Toward a proof of (\ref{comp1}) we shall use the divergence theorem on the set $\{u(x)>s\}$.
Choosing $\bfn= \nabla u/ |\nabla u|$ as inner unit normal vector, by (\ref{A9}) we have
\bel{duf}
\int_{\Sigma(s)} |\nabla u |\, d\ell~=~
\int_{\Sigma(s)} \nabla u \cdot \bfn \, d\ell~=~- \int_{\{u(x)>s\}} \Delta u\, dx~=~
 \int_{\{u(x)>s\}} \bigl( f(u) + c u_{x_1} - z\bigr)\, dx
\eeq
The inequality (\ref{comp1}) trivially holds when $s=1$, because in this case both sides vanish.    Let us decrease $s$ and check at what rate the two integrals increase.  
Taking the average values, and applying Jensen's inequality
for the convex function $y\mapsto 1/y$, for   any $c\in\R$ we obtain
\bel{dp1c}\bega{rl} \ds -{d\over ds} \psi(s)&\ds=~ f(s) \cdot \int_{\Sigma(s)} {1\over |\nabla u(x)|} \, d\ell(s)\\[4mm]
&\ds=~ f(s) \cdot \meas\bigl(\Sigma(s)\bigr) \cdot \avint_{\Sigma(s)} {1\over |\nabla u(x)|} \, d\ell(s)\\[4mm]
&\ds\geq~ f(s) \cdot \meas\bigl(\Sigma(s)\bigr) \cdot \left(\avint_{\Sigma(s)} |\nabla u(x)| \, d\ell(s)\right)^{-1}	\\[4mm]
&\ds=~ f(s) \cdot \Big[\meas\bigl(\Sigma(s)\bigr)\Big]^2 \cdot\left( \int_{\{u(x)>s\}} (f(u)+cu_{x_1}-z)\, dx \right)^{-1} \\[4mm]
&\geq ~ f(s) \cdot \Big[\meas\bigl(\Sigma(s)\bigr)\Big]^2 \cdot\bigl[\psi(s)+c(1-s)\bigr]^{-1},
\enda
\eeq
\bel{dp2c}\bega{rl} \ds -{d\over ds} \Psi(s)&\ds=~{f(s)\over U'(x(s))} ~=
~f(s) \cdot \left[ \int_{x(s)}^{+\infty} f(U(x)) \, dx +c(1-s)\right]^{-1}\\[4mm]
&\ds=~f(s) \cdot \bigl[\Psi(s)+c(1-s)\bigr]^{-1}.
\enda\eeq

Comparing (\ref{dp1c}) with (\ref{dp2c}), we see that
\begi
\item either~  $\psi(s)\geq \Psi(s)$,
\item or else ~
$\ds -{d\over ds} \psi(s)~\geq~-{d\over ds} \Psi(s)$.
\endi
Since $\psi(1)=\Psi(1)=0$, letting $s$ decrease from 1 to 0 by a comparison argument 
we conclude that  (\ref{comp1}) holds. \endproof

\begin{figure}[ht]
\centerline{\hbox{\includegraphics[width=15cm]{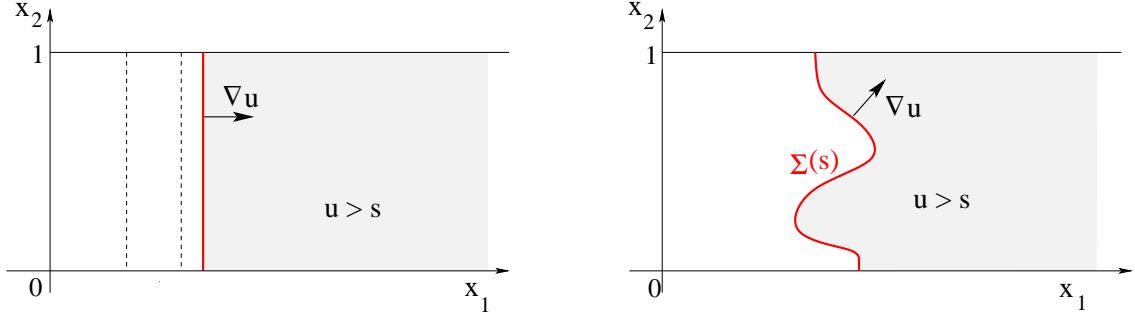}}}
\caption{\small  Left: the  regions where $u(x_1, x_2) = U(x_1)$ is $>s$. Right:
the region where
a general solution $u$ of (\ref{2dtv}) is $>s$. }
\label{f:co13}
\end{figure}

\section{The two-dimensional sharp interface limit}
\label{sec:7}
\setcounter{equation}{0}

We now return to the optimization problem introduced in Section~\ref{sec:1}, but 
with a possibly measure-valued dissipative source:
\bel{A5}
u_t~=~f(u) +\Delta u - \mu\,.\eeq
We are interested in the sharp interface limit, obtained by letting $\ve\to 0$ 
in the equation
\bel{A6} u_t~=~{1\over \ve}  f(u) + \ve\, \Delta u - \mu\,.\eeq
Notice that (\ref{A6}) can be derived from (\ref{A5}) 
simply by a rescaling of the independent variables
$t\mapsto \ve t$, $x\mapsto \ve x$. 
Here  $\mu$ as a (possibly measure-valued) non-negative control. 
For $\ve\approx 0$ we expect that the
solution to (\ref{A6}) will be a function taking values close to either 0 or 1 over most of its domain.  We thus seek to replace the controlled  parabolic equation (\ref{A6})
with  a control problem for a moving set.

The following notation will be used.  On the  
unit circumference  $S=\{ \xi\in \R^2\,;~|\xi|=1\}$ we use the arc-length measure, 
normalized so that $\int_S d\xi=1$.  For any vector $v=(v_1, v_2)\in \R^2$, the perpendicular 
vector (rotated by $90^o$) is 
$v^\perp = (- v_2, v_1)$.   By   ${\bf 1}_V$ we denote the characteristic function of a set $V\subset\R^2$, 
while $m_2(V)$ denotes its 2-dimensional Lebesgue measure.

\begin{figure}[ht]
\centerline{\hbox{\includegraphics[width=10cm]
{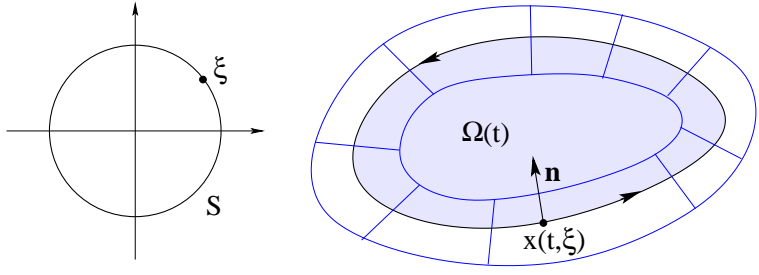}}}
\caption{\small At time $t$, the boundary $\partial \Omega(t)$ of the moving set 
is parameterized by the variable $\xi\in S^1$. An annulus of radius $\ve^{1/2}$
around this boundary is parameterized by 
$(\xi, s)\mapsto x(t,\xi) + s\bfn(t,\xi)$.  
}
\label{f:co29}
\end{figure}

Consider a set valued map $t\mapsto \Omega(t)\subset\R^2$.
For $t\in [0,T]$, let 
\bel{pom} \partial \Omega(t)~=~\{x(t, \xi)\,;~\xi\in S\}\eeq
be a $\C^1$ parameterization of the boundary  of $\Omega(t)$, oriented counterclockwise
(see Fig.~\ref{f:co29}).
We shall always assume that 
$x_\xi(t,x)\not= 0$ for all $t,\xi$, so that the unit inward normal vector  to $\Omega(t)$ at 
$x(t,\xi)$ is well defined by the formula
\bel{un}
\bfn(t,\xi)~\doteq~{x_\xi^\perp(t,\xi)\over |x_\xi(t,\xi)|}\,.\eeq
The normal velocity of the set boundary is given by the inner product
\bel{SE}\beta(t,\xi)~\doteq~\la \bfn(t,\xi)\,,\, x_t(t,\xi)\ra.\eeq

Throughout this section, we assume
that the source function $f$ satisfies {\bf (A2)} and {\bf (A4)}, so that Theorem~\ref{t:32} applies. 
In connection with the optimization problem {\bf (P2)} for a traveling wave, 
for every speed $c\geq c^*$ let $E(c)$ in (\ref{minc})
be  the minimum cost (\ref{J1}), among all 
measure-valued controls yielding a traveling profile with speed $c$.
One can extend $E$ to all values $c\in\R$ by setting 
$$E(c)\,=\, 0\qquad\hbox{for}\qquad c\leq c^*.$$
Integrating this cost along the boundary of a moving set, this leads to

\begin{definition} Consider a moving set $\Omega(t)$, with boundary parameterized as in 
(\ref{pom}). At each time $t\in [0,T]$,  the {\em instantaneous
effort} to achieve this motion is defined as
\bel{Eff}
\E(t)~\doteq
~ \int_{S}  E(\beta(t,\xi))\,\bigl|x_\xi(t,\xi)\bigr|\, d\xi\,.\eeq
\end{definition}

Given a convex function $\phi:\R_+\,\mapsto \R_+$ and two constants
$\kappa_1,\kappa_2\geq 0$, together with the optimization problem
{\bf (OP1)} introduced in Section~\ref{sec:1},  we now consider a problem of optimal control for the moving set $\Omega(t)$, $t\in [0,T]$. 
\begi
\item[{\bf (OP2)}] {\it Given an initial set $\Omega(0)=\Omega_0$, determine a controlled evolution
$t\mapsto \Omega(t)$ so that the total cost
\bel{cost1} J~\doteq~\int_0^T\phi\bigl(\E(t)\bigr)\, dt + \kappa_1\int_0^T m_2\bigl(\Omega(t)
\bigr)\, dt + \kappa_2\,m_2\bigl(\Omega(T)\bigr)\eeq
is minimized.}
\endi
A rigorous derivation of {\bf (OP2)} would require a study of the $\Gamma$-limit 
 of the functionals
\bel{Fen}\F_\ve(u)~\doteq~\int_0^T \int {[\ve \Delta u +\ve^{-1} f(u)-  u_t]_+\over u}\, dx\, dt
\eeq
as $\ve\to 0$.  Here $[s]_+= \max\{s,0\}$.  However,
this analysis is outside the scope of the present paper.  Here we only take some partial 
steps in this direction. The main result of this section shows that the cost $J$ 
at (\ref{cost1}) can be achieved as the limit of the cost (\ref{rtcost}), for a family
of solutions to the rescaled parabolic equations
 \bel{Aep} u^\ve_t~=~{1\over \ve}  f(u^\ve) + \ve\, \Delta u^\ve -u^\ve \alpha^\ve \,,
\qquad\qquad t\in [0,T], ~x\in\R^2.\eeq

\begin{theorem} \label{t:71}  Let $f$ satisfy the assumptions {\bf (A2)-(A3)}.
For $t\in [0,T]$, let $t\mapsto \Omega(t)\subset\R^2$ denote a moving set, whose boundary
admits a $\C^1$ parameterization as in (\ref{pom}).   
Moreover, assume that the normal velocity in (\ref{SE}) satisfies $\beta(t,\xi)\geq c^*$ for all
$t,\xi$.
Then there exists a family of control functions $\alpha^\ve$ and solutions
$u^\ve$ to (\ref{Aep})
such that the following two limits hold, uniformly for $t\in [0,T]$.
\bel{lep1}
\lim_{\ve\to 0} \Big\|u^\ve(t,\cdot) - {\bf 1}_{\Omega(t)} \Big\|_{\L^1}~=~0,\eeq
\bel{lep2}
\lim_{\ve\to 0}  \int_{\R^2}\alpha^\ve(t,x)\, dx~=~\int_{S^1}  E\bigl(\beta(t,\xi)\bigr)\,\bigl|x_\xi(t,\xi)\bigr|\, d\xi.\eeq

\end{theorem}  

{\bf Proof.}
{\bf 1.} By an approximation argument,  we can assume that 
the function $x=x(t,\xi)$ is smooth, and that  the normal speeds satisfy the strict inequality  
$\beta(t,\xi)> c^*$.   More precisely, we can choose constants $c_1, c_2, c_3$ 
such that 
\bel{c123} c^*\,<c_1\,<c_2\,<\, c_3, \qquad\qquad 
\beta(t,\xi)\,\in \,[c_2, c_3]\quad\forall  t\in [0,T],~
\xi\in S.\eeq
The solutions $u^\ve$ will be obtained by constructing suitable  lower and upper solutions
$u^\ve_- \leq u^\ve_+$.
\v
\begin{figure}[ht]
\centerline{\hbox{\includegraphics[width=15cm]
{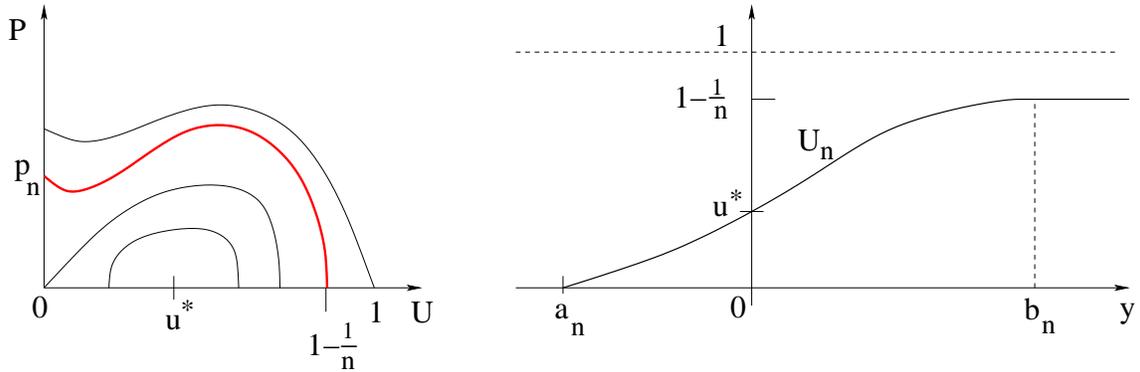}}}
\caption{\small The construction of lower solutions.  Left:
the trajectory of (\ref{UPc1}) through the point $(1-n^{-1},0)$ crosses the 
$P$-axis at some point $(0, p_n)$.  Right: the corresponding  traveling profile at (\ref{Un1})-(\ref{Un2})
is uniquely determined by requiring that $U_n(0)= u^*$. 
By construction we have  $U_n'(a_n) = p_n>0$, while $U_n'(b_n)=0$.}
\label{f:co20}
\end{figure}

{\bf 2.} Toward the construction  of lower solutions (see Fig.~\ref{f:co20}),
let $c_1>c^*$ be given. For every
$n\geq 1$ large enough, the trajectory of the system
\bel{UPc1}\left\{ \bega{rl} U'&=~P,\\[1mm]
P'&=~- f(U) - c_1 P,\enda\right.\eeq
that goes through the point $(1-n^{-1},\,0)$ will cross the $P$-axis at a point $(0, p_n)$,
with 
\bel{ann} p_n\,>\,0,\qquad \lim_{n\to \infty} p_n ~=~ \bar p~ >~0.\eeq
This yields a traveling profile $U_n= U_n(y)$ which satisfies
\bel{Un1} U_n(0)= u^*,\qquad \qquad U_n'' + c_3 U_n' +f(U_n)~=~0\quad \hbox{for} ~~y\in [a_n, b_n],  \eeq
\bel{Un2} \left\{ \bega{rl} U_n(a_n)&= ~0,\\[2mm]
U'_n(a_n)&= ~p_n~>~0,\enda\right.\qquad \left\{ \bega{rl} 
U_n(b_n)&= ~1-{1\over n},\\[2mm]
U'_n(b_n)&= ~0,\enda\right.\eeq
for some values $a_n<0<b_n$.
We can extend it outside the interval $[a_n, b_n]$ by setting
\bel{Un3} U_n(y)~=~\left\{ \bega{cl} 0\quad &\hbox{if}\quad y\leq a_n\,,\\[2mm]
 1-{1\over n}\quad &\hbox{if}\quad y\geq b_n\,.\enda\right.\eeq
\v

\begin{figure}[ht]
\centerline{\hbox{\includegraphics[width=15cm]
{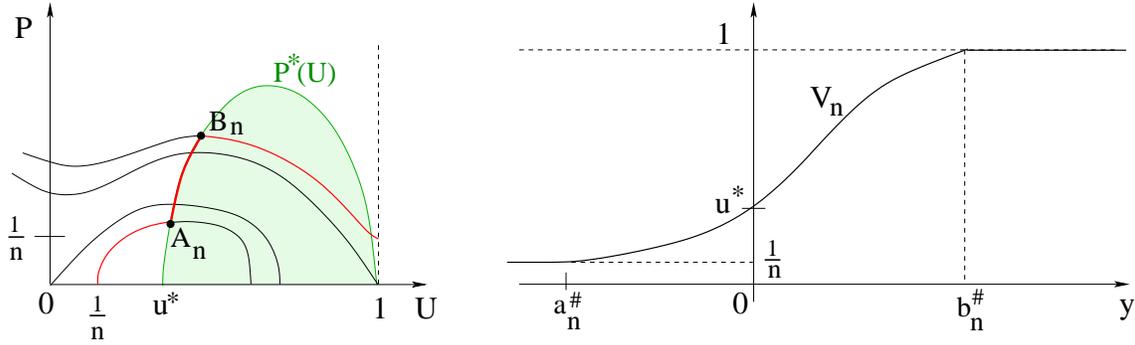}}}
\caption{\small The construction of upper solutions.  Left:
a concatenation of a trajectory  of (\ref{T2}) through $(n^{-1},0)$, followed by an arc along the curve where
$P=P^*(U)$, followed by  a trajectory  of (\ref{T2}) through $(1, n^{-1})$.
 Right: the corresponding  traveling profile in (\ref{Vn1})-(\ref{Vn2}), shifted so that $V_n(0)= u^*$. 
By construction we have  $V_n'(a'_n) = 0$, while $V_n'(b'_n)>0$.}
\label{f:co21}
\end{figure}

{\bf 3.} Toward the construction of upper solutions,
consider again the 1-dimensional equation
\bel{cwe2}u_t~=~f(u)  + u_{xx}- u\, \alpha \,,\eeq
where $\alpha=\alpha(t,x)\geq 0$ is the control function. 
For a given wave speed $c> c^*$, a traveling profile $u(t,x)= U(x-ct)$ is an upper solution
provided that
\bel{uus}U'' +c U' + f(U) - \alpha\, U~\leq~0.\eeq 
Assuming that $f$ satisfies {\bf (A2)-(A3)}, we 
consider the path $\gamma_n$ obtained by concatenating the following three curves
(see Fig.~\ref{f:co21}, left)
\begi\item   The trajectory of (\ref{T2}) starting at $( n^{-1},0)$, up to the point $A_n$
where
it intersects the curve $P=P^*(U)\doteq \sqrt{Uf(U)}$.
\item The trajectory of (\ref{T2}) ending at $(1, n^{-1})$,  continued backward 
up to a point $B_n$ along the curve where $P=P^*(U) $.  
\item The portion of the curve $P=P^*(U)$  between $A_n$ and $B_n$.
\endi
As shown in Fig.~\ref{f:co21}, right,
this yields a traveling profile for (\ref{cwe2}), say $u(t,x) = V_n(x-ct)$, 
 which satisfies
\bel{Vn1} V_n(0)= u^*,\qquad \qquad V_n'' + c V_n' +f(V_n)+V_n \alpha_n~=~0\quad \hbox{for} ~~y\in [a^\sharp_n, b^\sharp_n],  \eeq
\bel{Vn2} \left\{ \bega{rl} V_n(a^\sharp_n)&= ~{1\over n},\\[2mm]
V'_n(a^\sharp_n)&= ~0,\enda\right.\qquad \left\{ \bega{rl} 
V_n(b^\sharp_n)&= ~1,\\[2mm]
V'_n(b^\sharp_n)&> ~0,\enda\right.\eeq
for some values $a^\sharp_n<0<b^\sharp_n$ and a suitable control $\alpha_n\geq 0$.
We can extend $V_n$  outside the interval $[a^\sharp_n, b^\sharp_n]$ by setting
\bel{Vn3} V_n(y)~=~\left\{ \bega{cl} {1\over n}\quad &\hbox{if}\quad y\leq a^\sharp_n\,,\\[2mm]
 1\quad &\hbox{if}\quad y\geq b^\sharp_n\,.\enda\right.\eeq
 We remark that the above profiles $V_n$, as well as the interval $[a^\sharp_n, b^\sharp_n]$,
 all depend on the wave speed $c$.   To be reminded of this fact, we shall use the notations
 $V_n^{c}$, $a^\sharp_n(c)$, $b^\sharp_n(c)$.  
 
 In connection with (\ref{c123}) we observe that, as long as the speed
 $c\in [c_2, c_3]$ remains in a bounded interval, also the intervals $\bigl[a^\sharp_n(c), 
 b^\sharp_n(c)\bigr]$ remain uniformly bounded
\v
{\bf 4.} Using the above traveling profiles $U_n, V_n$, we are now ready to construct 
sequences of upper and lower solutions.

Choose a sequence $\ve_n\downarrow 0$ such that, for every  $n\geq 1$,
\bel{venc}{1\over \sqrt{\ve_n}} ~\geq~(b_n - a_n) + \bigl(b^\sharp_n(c)- a^\sharp_n(c)\bigr)
\quad\qquad\forall ~ c\in [c_2, c_3].\eeq
Define the rescaled profiles
\bel{TUn}
\Tilde U_n(y)~=~U_n\left( y- \sqrt{\ve_n}\over \ve_n\right),\qquad\qquad \Tilde V_n^c(y)
~\doteq~V_n^{c}\left( y+ \sqrt{\ve_n}\over \ve_n\right).\eeq
Notice that, by the definitions of $U_n$ and $V_n^{(c)}$, this implies
\bel{TUVn}
\Tilde U_n(y)~=~\left\{ \bega{cl} 0\quad &\hbox{if}~~y\leq 0,\\[1mm]
1-n^{-1}\quad &\hbox{if}~~y\geq  2\sqrt{\ve_n}\,,\enda\right.
\qquad\qquad \Tilde V^{c}_n(y)~=~\left\{ \bega{cl} n^{-1} \quad &\hbox{if}~~y\leq - 2\sqrt{\ve_n},\\[1mm]
1\quad &\hbox{if}~~y\geq  0\,.\enda\right.\eeq

\begin{figure}[ht]
\centerline{\hbox{\includegraphics[width=15cm]
{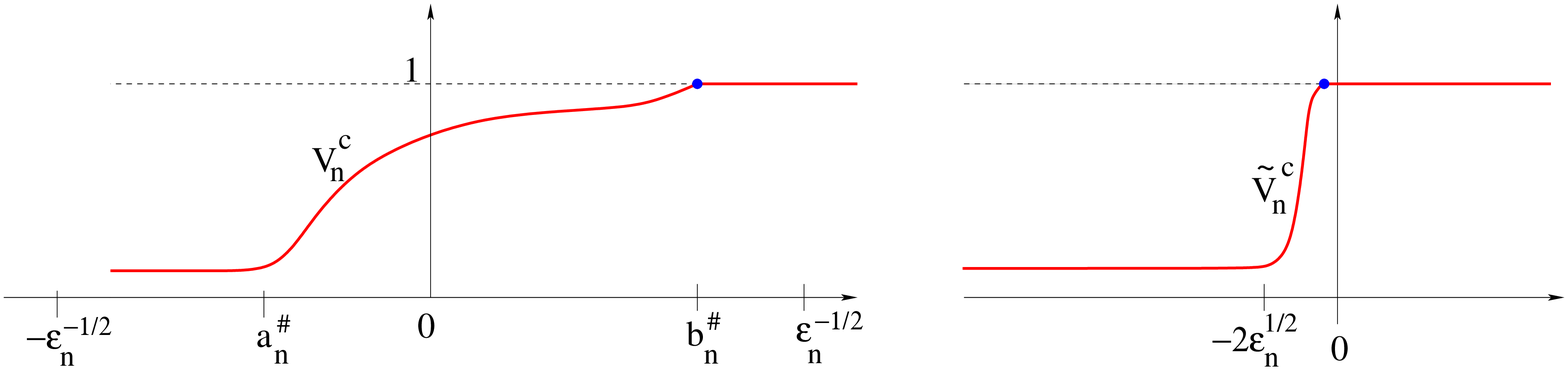}}}
\caption{\small The traveling profile $V^c_n$, providing an upper solution, and its rescaled version $\Tilde V^c_n$ at (\ref{TUVn}).}
\label{f:sc18}
\end{figure}

Recalling the parameterization (\ref{pom}) of the boundary $\partial\Omega(t)$, 
consider the annular domain
\bel{Dnt}
\D_n(t)~\doteq~\Big\{ x(t,\xi)+ y\,\bfn(t,\xi)\,;\quad \xi\in S, \quad |y|\leq 2\sqrt{\ve_n}\Big\}.
\eeq
Thanks to our earlier  assumptions on the map $(t,\xi)\mapsto x(t,\xi)$,
for all $\ve_n>0$ small enough the map 
$$(\xi,y)~\mapsto~x(t,\xi)+ y\,\bfn(t,\xi)\qquad\qquad x(t,\xi)+ y\,\bfn(t,\xi)$$
has a smooth inverse, for every $t\in [0,T]$.

We now define a lower solution $u^-_n$   by setting
\bel{ulow}u^-_n\Big(x(t,\xi)+ y\,\bfn(t,\xi)\Big)~\doteq~\Tilde U_n(y)\qquad\qquad 
\forall ~|y|\leq 2\sqrt{\ve_n},\eeq
and extending $u_n^-$ outside the annulus $\D_n(t)$  as a constant function.
Namely: $u_n^-(t,x)=0$
 in the interior of 
$\Omega(t)$, while  $u_n^-(t,x)=1-n^{-1}$ outside $\Omega(t)$.
This is possible because of (\ref{TUVn}).

Similarly, we define an upper solution $u^+_n$
by setting
\bel{uupp}u^+_n\Big(x(t,\xi)+ y\,\bfn(t,\xi)\Big)~\doteq~\Tilde V^{\beta(t,\xi)}_n(y),\eeq
and extending $u_n^-$ outside the annulus $\D_n(t)$  as a constant function.
Namely: $u_n^+(t,x)=n^{-1}$
 in the interior of 
$\Omega(t)$, while  $u_n^-(t,x)=1$ outside $\Omega(t)$. 
Again, this is possible because of (\ref{TUVn}).
%

\v
By construction, it is clear that $0\leq u^-_n(t,x)\leq u^+_n(t,x)\leq 1.$
Moreover, $$u_n^-(t,\cdot)\to {\bf 1}_{\Omega(t)}
\qquad\hbox{in}\qquad \L^1(\R^2).$$
However, we only have 
 $$u_n^+(t,\cdot)\to {\bf 1}_{\Omega(t)}
\qquad\hbox{in}\qquad \L^1_{loc}(\R^2),$$
because $u_n^+(t,x)\geq n^{-1}$ and hence this upper solution is not integrable on $\R^2$.

To cope with this problem, observing that the minimum between two upper solutions is
an upper solution, we can proceed as follows.
Let $R>0$ be a radius large enough so that  all sets $\Omega(t)$ remain inside  
the disc $B(0,R)\subset\R^2$.  
Consider the radially symmetric function 
\bel{phi}\phi(x)~\doteq~\begin{cases}1 \quad &|x|\leq R,\\
e^{-|x|+R}  \quad &|x|\geq R \end{cases}\eeq
Notice that, for $\ve$ small enough, this is a time-independent upper solution of
$$u_t~=~\ve \Delta u +{1\over\ve}  f(u),$$
integrable over the entire plane $\R^2$.

Replacing the functions $u_n^+$ with
$$v_n(t,x) ~=~\min ~\Big\{u_n^+(t,x) ,\, \phi(x)\Big\},$$
we obtain  a new sequence of upper solutions. We claim that this sequence
converges to ${\bf 1}_{\Omega(t)}$ in $\L^1(\R^2)$, for every $t\in [0,T]$.
Indeed, for every $n\geq 1$ sufficiently large one has
\bel{conv} \bega{rl}\ds\int_{\R^2}|v_n(t,x)-{\bf 1}_{\Omega(t)}|\,dx~=~\int_{B(0,R+\ln n)\setminus \Omega(t) }u^+_n(t,x)\, dx ~ + ~ \int_{|x|>R+\ln n}  \phi(x)\, dx \\ [4mm]
\ds ~ \leq ~\meas (\D_n(t)\setminus\Omega(t))+ {\pi(R+\ln n)^2\over n} +
 \int_{|x|>R+\ln n} e^{-|x|+R}\, dx,
 \enda\eeq
and each term on the right hand side of (\ref{conv}) goes to zero as $n \to \infty$.
\v
{\bf 5.} For each $n\geq 1$, we now consider the cost of a control 
which renders $v_n$ an upper solution.
The smallest control function $\alpha$ that fulfils this requirement is
\bel{anmin}
\alpha_n~\doteq ~{\bigl[\ve_n \Delta v_n +  \ve_n^{-1} f(v_n) - (v_n)_t\bigr]_+\over 
v_n}\eeq
By construction, we already know that $v_n$ satisfies
$$\ve_n\, \Delta v_n + \ve_n^{-1} f(v_n) ~\leq~0~=~(v_n)_t\qquad \hbox{for}~~x\notin \D_n(t).$$
It thus remains to estimate the integral
\bel{intan}\int_{\D_n(t)}
\alpha_n(t,x)\, dx,\eeq
and show that it converges to the right hand side of (\ref{lep2}).
\v
%
{\bf 6.} 
Over the set $\D_n(t)$, we shall use the coordinates $(\xi,y)\in S\times [-2\sqrt{\ve_n},
\, 2+\sqrt{\ve_n}]$ corresponding to the point
\bel{cvar}x~=~x(t,\xi) + y \,\bfn(t,\xi)\,.\eeq

Computing the Laplacian of $v_n$ in terms of the coordinates
$\xi,y$, and calling $r= r(t,\xi,y)$ is the local radius of curvature (which is uniformly positive throughout the domain), we find
\bel{Devn}\bega{rl}
\Delta_x v_n&\ds=~(v_n)_{yy} + {1\over r} (v_n)_y + {(v_n)_{\xi\xi} - x_{\xi\xi} \cdot (v_n)_\xi\ 
	\over |x_\xi|^2}\\[3mm]
&\ds=~{1\over\ve_n^2} \bigl(V_n^{\beta(t,\xi)}\bigr)'' \left(y+\sqrt{\ve_n}\over \ve_n\right) +
\O(1)\cdot {1\over \ve_n} + \O(1)\,.
\enda\eeq
On the other hand,
\bel{Dtvn}
\bega{rl}
\partial_t v_n&\ds=~- {\beta(t,\xi)\over \ve_n} \bigl(V_n^{\beta(t,\xi)}\bigr)' \left(y+\sqrt{\ve_n}\over \ve_n\right) +
\O(1).
\enda\eeq
Using (\ref{Vn1}), with the speed $c=\beta(t,\xi)$, one obtains
$$\bigl(V_n^{\beta(t,\xi)}\bigr)'' + \beta(t,\xi)\bigl(V_n^{\beta(t,\xi)}\bigr)' + f\bigl(V_n^{\beta(t,\xi)}\bigr) + V_n^{\beta(t,\xi)}\cdot \alpha_n~=~0.$$
Combining the above estimates, we obtain
\bel{71}\ve_n \Delta_x v_n + {1\over \ve} f(v_n) - \partial_t v_n~=~{1\over \ve_n}  V_n^{\beta(t,\xi)}\cdot \alpha_n + \O(1),
\eeq
where $\alpha_n$ is the optimal control on the portion of curve from $A_n$ to $B_n$
in Fig.~\ref{f:co21}.
\v
{\bf 7.} We now integrate (\ref{71}) over the entire domain $\D_n$. Computing the Jacobian determinant of the transformation
(\ref{cvar}), we obtain
$$dx_1 dx_2~=~\bigl| x_\xi(\xi,y)\bigr| \cdot  \ve_n \, d\xi d y. $$
Therefore,  at any given time $t\in [0,T]$, there holds
\bel{aln}\bega{rl}\ds \int_{\R^2} \alpha_n\, dx&\ds=~\int_{\D_n} {\bigl[\ve_n \Delta v_n +  \ve_n^{-1} f(v_n) - (v_n)_t\bigr]_+\over 
	v_n}
\, dx\\[4mm]
&\ds =~\int_{S\times [-2\sqrt {\ve_n},~2\sqrt{\ve_n}]} \left[ {\alpha_n(\xi,y)
	\over\ve_n} + \O(1)\right]\, \bigl| x_\xi(\xi,y)\bigr|  \ve_n \,d\xi d y\\[4mm]
&\ds =~\int_S E\bigl(\beta(t,\xi)\bigr) |x_\xi(t,\xi)\bigr|\, d\xi + \O(1)\cdot \ve_n\,.
\enda \eeq
Taking the limit 
as $\ve_n\to 0$, this yields (\ref{lep2}).
\v
{\bf 8.} Having constructed  a sequence of lower solutions $u_n^-$ and of upper solutions
$v_n$ which both converge to the characteristic function ${\bf 1}_{\Omega(t)}$, by 
a comparison argument we obtain a sequence of solutions to
\bel{unsol}  u_{n,t}~=~{1\over \ve_n}  f(u_n) + \ve_n\, \Delta u_n -u_n \alpha_n \,,
\qquad\qquad t\in [0,T], ~x\in\R^2, \eeq
with $u_n^-\leq u_n \leq v_n$. By (\ref{aln}), these solutions satisfy
\bel{lu1}
\lim_{n\to\infty} \Big\|u_n(t,\cdot) - {\bf 1}_{\Omega(t)} \Big\|_{\L^1}\,=~0,\qquad \quad
\lim_{n\to\infty}  \int_{\R^2}\alpha_n(t,x)\, dx~=~\int_{S}  E(\beta(t,\xi))\,\bigl|x_\xi(t,\xi)\bigr|\, d\xi.\eeq
This achieves the proof.
\endproof

\begin{remark} {\rm We expect that an entirely similar result could be proved in the 
case where 
$g(u,\alpha)= \alpha$, and $f$ satisfies either (\ref{f1}) or (\ref{f2}).   In this case, the optimal traveling profiles are the ones described in 
Theorem~\ref{t:31}, while the formulas (\ref{Aep}) and (\ref{Fen}) 
should be replaced respectively by
$$ u^\ve_t~=~{1\over \ve}  f(u^\ve) + \ve\, \Delta u^\ve -\alpha^\ve \,,
\qquad\qquad \F_\ve(u)~\doteq~\int_0^T \int [\ve \Delta u +\ve^{-1} f(u)-  u_t]_+\, dx\, dt \,.
$$
}
\end{remark} 
\v

\end{document}